\documentclass[12pt]{amsart}
\usepackage{epsf, amssymb, xy, latexsym, amsmath, graphics}
\newtheorem{theorem}{Theorem}
\newtheorem{corollary}{Corollary}
\newtheorem{lemma}{Lemma}
\newtheorem{proposition}{Proposition}
\newcommand{\QED}{{ \hfill \qed}} 
\newcommand{\comp}{\, {{\rm o}}\,} 
\def \1{^{-1}}
\def \2{^{-2}}
\def \3{^{-3}}
\def \pr{^\prime}
\def \prpr{^{\prime\prime}}
\renewcommand{\phi}{\varphi}

\def \b{\beta}
\def \c{\gamma}

\def \q{\kappa}
\def \t{\tau}
\def \s{\sigma}
\newcommand{\ca}{{\mathcal A}}
\newcommand{\cb}{{\mathcal B}}
\newcommand{\cc}{{\mathcal C}}
\newcommand{\cd}{{\mathcal D}}
\newcommand{\ce}{{\mathcal E}}

\newcommand{\cm}{{\mathcal M}}

\newcommand{\cs}{{\mathcal S}}
\newcommand{\co}{{\mathcal O}}
\newcommand{\cq}{{\mathcal Q}}
\newcommand{\C}{{\mathbb  C}} 
\newcommand{\N}{{\mathbb  N}} 
\newcommand{\Z}{{\mathbb  Z}}

\newcommand{\pline}{{{\mathbb P}^1}} 
\newcommand{\proj}{{{\mathbb  P}^2}}  
\newcommand{\ball}{{\mathbf B_{\mathbf 2}}} 
\newcommand{\oball}{{\mathbf B_{\mathbf 1}}} 
\newcommand{\pball}{\oball\!\times \oball} 
\newcommand{\pp}{\pline\times\pline} 
\newcommand{\fg}[1]{{\pi_1(\proj\moins #1)}} 
\newcommand{\orbfg}{\pi_1^{orb}} 
\newcommand{\ok}{\to}

\newcommand{\moins}{\backslash}
\newcommand{\sol}{\langle} 
\newcommand{\sag}{\rangle} 
\newcommand{\ssag}{{\rangle\!\rangle}} 
\newcommand{\ssol}{\langle\!\langle} 
\begin{document}
\title{Covering Relations Between Ball-Quotient Orbifolds} 
\author{A.Muhammed Uluda{\u g}}
\date{}
\maketitle
\begin{abstract} 
Some ball-quotient orbifolds are related by covering maps.
We exploit these coverings to find infinitely many 
orbifolds on $\proj$ uniformized by the complex 2-ball $\ball$ and 
some orbifolds over K3 surfaces  uniformized by $\ball$.
We also give, along with infinitely many reducible examples,
an infinite series of irreducible curves along which $\proj$ is 
uniformized by the product of 1-balls $\pball$.
\end{abstract}
\section{Introduction}
\label{intro}
Let $M$ be a smooth algebraic surface of general type.
By the Hirzebruch proportionality theorem~\cite{hirzebruch0}
if $M$ is covered by the complex 2-ball $\ball$
then its Chern numbers satisfies the equality $c_1^2(M)=3e(M)$.
Conversely if the Chern numbers of $M$ satisfy this equality 
then the universal covering of $M$ is $\ball$ by Yau's theorem~\cite{yau}.
\par
One way to discover surfaces of general type $M$ with $c_1^2(M)=3e(M)$
is to construct them as  finite branched Galois coverings $\phi:M\ok X$ where $X$ 
is the blow-up of $\proj$ at some points.
This approach was first used  by 
Hirzebruch~\cite{hirzebruch},~\cite{hirzebruch2},
developed further in~\cite{hofer},~\cite{barthel} and more recently 
in~\cite{holzapfel3},\cite{holzapfel}. 
The corresponding lattices acting on $\ball$ turned out to be commensurable 
with the ones obtained from the study of the hypergeometric differential 
equations (see~\cite{deligne2} and ~\cite{yoshida2}). 
\par
A branched Galois covering $\phi:M\ok X$ endows $X$ with a map $\beta_\phi:X\ok \N$ 
sending $p\in X$ to the order of the isotropy group above $p$. 
The pair $(X,\beta_\phi)$ is an orbifold, {and $M$ is   
a \textit{uniformization} of $(X,\beta_\phi)$. In case the degree of $\phi$ is finite, $M$ 
is also called  a \textit{finite uniformization}, and if $M$ is simply connected,
it is called the \textit{universal uniformization} of  $(X,\beta_\phi)$. 
Alternatively, when finiteness or universality can be understood from the context,
one says: \textit{$(X,\beta_\phi)$ is uniformized by M}.}
\par
For an orbifold $(X,\beta)$ it is possible to define 
orbifold Chern numbers
$c_1^2(X,\beta)$ and $e(X,\beta)$ in such a way that for any finite uniformization 
$\phi:M\ok X$ with $\beta_\phi=\beta$ one has 
$$
c_1^2(M)=\mbox{deg}(\phi) c_1^2(X,\beta),\qquad  e(M)=\mbox{deg}(\phi) e(X,\beta)
$$  
On the other hand, by results of \cite{kobayashi} an orbifold of general type  $(X,\beta)$ is covered by   
$\ball$ if the equality 
\begin{equation}
\label{cchern}
c_1^2(X,\beta)=3e(X,\beta)
\end{equation} 
is satisfied. Hence, in order to discover lattices acting on $\ball$ 
it suffices to find orbifolds satisfying (\ref{cchern}).
In this article, this idea is applied to rediscover some orbifolds discovered in~\cite{holzapfel}. 
A study of some coverings (in the orbifold sense) of these orbifolds gave the following result.
\begin{theorem}
\label{mr}
There exists infinite series of  pairwise non-isomorphic 
orbifolds on $\proj$, uniformized by the 2-ball $\ball$.
Corresponding lattices in ${\rm SU}(1,2)$ are all arithmetic.
\end{theorem}
In contrast with many known examples of orbifolds $(\proj,\beta)$ 
uniformized by $\ball$,
to the author's knowledge there is just one orbifold $(\proj, \beta)$ 
(alleged to Hirzebruch in~\cite{yoshida}) 
that is known to be uniformized by $\pball$. 
In Section \ref{sectioncoverings} below, infinitely many examples of such orbifolds are given.
In particular, the following result is proved.
\begin{theorem}\label{thm2}
For $m>0$ odd, let $Q_m$ be the irreducible  curve given by the equation
$x^{m/2}+y^{m/2}+z^{m/2}=0$.
Define $\b_m:\proj\ok\N$ by $\b_m(x)=2$ if 
$x\in Q_m\moins\mbox{\rm sing}(Q_m)$, $\b_m(x)=2m$ if 
$x\in\mbox{\rm sing}(Q_m)$ and $\b_m(x)=1$ otherwise. 
Then the orbifold $(\proj,\b_m)$ is uniformized by $Q_m\times Q_m$.
Hence, for $m>3$ the universal uniformization of $(\proj,\b_m)$ is $\pball$.
{The uniformization $Q_m\times Q_m\ok (\proj,\b_m)$ is of degree $2m^2$.} 
\end{theorem}
The curve $Q_1$ is a smooth quadric, and it is well known that $(\proj,\b_1)$ is uniformized 
by $\pp$. The curve $Q_3$ is the nine-cuspidal sextic, it was shown in \cite{kaneko2} that 
$\C^2$ is the universal uniformization of $(\proj, \b_3)$. 
{For a generalization of Theorem~\ref{thm2}
to higher dimensional projective spaces, see~\cite{uludag}.}
\section{Orbifolds}
We shall mostly follow the terminology settled in~\cite{yy}.
Let $M$ be a connected complex manifold, $G\subset \mbox{Aut}(M)$ a properly
discontinuous subgroup and put $X:=M/G$. 
Then the projection $\phi:M\ok X$ is a branched Galois 
covering  endowing $X$ with a map $\beta_\phi:X\ok\N$ defined by 
$\b_\phi(p):=|G_q|$ where 
$q$ is a point in $\phi^{-1}(p)$ and $G_q$ is the isotropy subgroup of $G$ at $q$. 
In this setting, the pair $(X,\b_\phi)$ is said to be uniformized by $\phi:M\ok (X,\b_\phi)$.
An \textit{orbifold} is a pair $(X,\b)$ of an irreducible normal analytic space $X$ with a 
function $\b:X\ok \N$ such that  the pair $(X,\b)$ is locally finitely uniformizable. 
Let $(X,\b)$ and $(X,\c)$ be two orbifolds with $\c|\b$, and let 
$\phi:(X\pr,1)\ok (X,\c)$ be a uniformization of $(X,\c)$, e.g. $\b_\phi=\c$. 
Then  $\phi:(X\pr,\b\pr)\ok (X,\b)$ is called an \textit{orbifold covering}, where 
$\b\pr:=\b\comp\phi/\c\comp\phi$. 
The orbifold $(X\pr,\b\pr)$ is called the 
\textit{lifting of $(X,\b)$ to the uniformization $X\pr$ of $(X,\c)$}.
\par
Let $(X,b)$ be an orbifold,  
$B_\b:={\rm supp}(\b-1)$ and let $B_1,\dots,B_n$ be the irreducible components of $B_\b$.
Then $\b$ is constant on $B_i\moins {\rm sing}(B_\b)$; so let $b_i$ be this number.
The \textit{orbifold fundamental group} $\orbfg(X,\b)$ of $(X,\b)$ is the group defined by
$\orbfg(X,\b):=\pi_1(X\moins B_\b)/\ssol\mu_1^{b_1},\dots, \mu_n^{b_n}\ssag$
where $\mu_i^{b_i}$ is a meridian of $B_i$ and $\ssol\ssag$ denotes the normal closure.
An orbifold $(X,\b)$ is said to be \textit{smooth} if $X$ is smooth. 
In case $(X,\b)$ is a smooth orbifold the map $\beta$ is determined by the numbers 
$b_i$; in fact  $\b(p)$ is the order of the local orbifold fundamental group at $p$. 
In case $\dim X=2$ the orbifold condition (i.e. locally finitely uniformizability) is equivalent to 
the finiteness of the local fundamental groups.
For example if $x=B_i\cap B_j$ is a simple node of $B_\b$  then $\b(p)=b_ib_j$.
For some other singularities of $B_\beta$ one has
\begin{lemma}\label{tangency}
Let $(X,\b)$ be an orbifold where $X$ is a smooth complex surface and $p\in X$. 
(i) If $p=B_i\cap B_j\cap B_k$ is a 
transversal intersection of smooth branches of $B_\beta$  then 
$b_i\1+b_j\1+b_k\1>1$ and $\b(p)=4[{b_i\1}+{b_j}\1+{b_k}\1-1]^{-2}$
(ii) If $p= B_i\cap B_j$ is a tacnode of
$B_\b$ then  $b_i\1+b_j \1>1/2$ and 
$\b(p)=2[{b_i\1}+{b_j\1}-{2\1}]^{-2}$
(iii)  If $p\in B_i$ is a simple cusp of $B_\b$ then  $b_i\1> 6\1$ and 
$\b(p)=\frac{2}{3}[{b_i\1}-{6 \1}]^{-2}$.
\end{lemma}

\begin{proof}
Part (i) is well known, see e.g.~\cite{yoshida}.
Now let $B_1,\, B_2,\, B_3\, \subset\C^2$ be respectively the lines $x=0$, $x=y$ and $y=0$.
By part (i) the pair $(\C^2,\b)$ is an orbifold where $\b:\C^2\ok \N$ is the function  
$$
\b(p)=\left\{\begin{array}{ll}
                b_i,                            & p\in B_i\moins\{(0,0)\} \\
                4[{b_1\1}+{b_2}\1+{b_3}\1-1]^{-2},      & p=(0,0)\\
                1,                              & {\rm otherwise}
        \end{array}\right.
$$
for integers $b_1\, ,b_2\, ,b_3$ satisfying  $b_1\1+b_2\1+b_3\1>1$. 
To prove (ii) put $b_3=2$ and consider the branched Galois covering 
$\phi: (x,y)\in\C^2\ok (x,y^2)\in\C^2$. 
One has $\b_\phi(p)=b_3=2$ for $p\in B_3$ and $\b_\phi(p)=1$ otherwise.
Let $B_1\pr, B_2\pr \subset \C^2$ be respectively the line $x=0$ and the curve $x=y^2$.
Then $\phi(B_1\pr)=B_1$ and $\phi(B_2\pr)=B_2$. Let $\b\pr:\C^2\ok \N$ be the function 
which takes the value $b_i$ on $B_i\pr\moins\{(0,0)\}$, 
the value $2[{b_i\1}+{b_j\1}-{2\1}]^{-2}$ on $(0,0)$ and the value $1$ otherwise. 
Then $\phi:(\C^2,\b\pr)\ok (\C^2,\b)$ is an orbifold covering, which proves (ii). 
To prove (iii) one applies the above argument with $\phi: (x,y)\in\C^2\ok (x^3, y^2)$, 
$b_1=3$ and $b_2=2$. (Setting $B_1\pr$ to be the curve $x^3=y^2$, one has 
$\phi(B_2\pr)=B_2$).  
\end{proof}

\medskip\noindent\textbf{Conventions.}
We shall almost exclusively be concerned with orbifolds 
$(X,\b)$ with $X$ being a smooth algebraic surface.
In most cases $X$ will be the projective plane $\proj$. 
Since in this case  $\b$ is determined by its values $b_i$ on $B_i\moins {\rm sing}(B_\b)$,
a smooth orbifold can alternatively be defined as a pair $(X,B)$ where  
$B:=b_1B_1+\dots b_nB_n$ is a divisor on $X$ with  $b_i\geq 1$. 
For an orbifold $(X,B)$ the corresponding map $X\ok\N$ will be denoted by $\b_B$.
{The \textit{locus} of the orbifold $(X,B)$ is the hypersurface $B\subset X$. } 

\subsection{Orbifold Chern numbers.}
Let $(\proj,B)$ be an orbifold where $B=b_1B_1+\dots+b_nB_n$ with  
$B_i$ being an irreducible curve of degree $d_i$.

\medskip\noindent\textbf{Definition.}
The \textit{orbifold Chern numbers of} $(\proj,B)$ are defined as 
$$
c_1^2(\proj,B):=\bigl[-3+\sum_{1\leq i\leq n}d_i\left(1-{b_i}\1\right)\bigr]^2
$$
$$
e(\proj,B):=3-\sum_{1\leq i\leq n}\left(1-{b_i}\1\right)e(B_i\moins{\rm sing}(B))
-\sum_{p\in{\rm sing}(B)}\left(1-{\b_B(p)}\1\right)
$$

If $(M,1)\ok (\proj, \b)$ is a finite uniformization  of degree $d$
then the Chern numbers of $M$ are given by 
$e(M)=de(\proj,B)$ and $c_1^2(M)=dc_1^2(\proj,B)$.
\par
{In dimension 2, an orbifold $(X,\b)$ is said to be of \textit{general type} 
if it is uniformized by a surface of general type. In the context of the following theorem,
this simply means that $(X,\b)$ is \textit{not} uniformized by $\proj$.} 

\begin{theorem}[Kobayashi, Nakamura, Sakai \cite{kobayashi}]\label{kns}
Let $(\proj,B)$ be an orbifold of general type. 
Then $c_1^2(\proj,B)\leq 3e(\proj,B)$, the equality holding if and only if 
$(\proj,B)$ is uniformized by $\ball$.
\end{theorem}
In fact, the KNS theorem is proved in greater generality then its version stated above; in particular it is
valid for orbifolds $(\proj,\b)$ with at worst ``log-canonical singularities'' and $\b$ being a function 
$\proj\ok \N\cup{\infty}$. This implies that in Lemma~\ref{tangency} (i) one may have
$b_i\1+b_j\1+b_k\1=1$ with $\b(p)=\infty$, in (ii) one may have $b_i\1+b_j \1=1/2$ with 
$\b(p)=\infty$ and in (iii) one may have $b_i\1=1/6$ with 
$\b(p)=\infty$.
\par
Let $M$ be an algebraic surface with $\pball$ as the universal covering. 
Then by the Hirzebruch proportionality~\cite{hirzebruch0} one has $c_1^2(M)=2e(M)$. 
Similarly if an orbifold $(X,\b)$ is uniformized by $\pball$ then by 
Selberg's theorem the corresponding transformation 
group has a torsion-free normal subgroup of finite
index $d$, which implies that $(X,\b)$ admits a finite uniformization $M\ok X$ of degree $d$. 
Since $M$ is uniformized by  $\pball$, one has  
$dc_1^2(X,\b)=c_1^2(M)=2e(M)=2de(X,\b)$.

\section{The Apollonius configuration}\label{apoconf}
Let $A_n:=Q\cup T_1\cup\dots \cup T_n$ 
be an arrangement  consisting of a smooth quadric $Q$ with $n$  
distinct tangent lines of $Q$. 
Since there are only two tangent lines to a quadric from a point 
$\in\proj\moins Q$, the lines $T_1,\dots ,T_n$ meets each other one by one. 
The configuration space of $A_n$'s is naturally identified
with the configuration
space $\cm_n$ of $n$ distinct points in $\pline$, via the contact 
points of the tangent lines with the quadric $Q\simeq \pline$.
Since the space $\cm_n$ is connected,
any two arrangements $A_n$ with $n$ fixed are isotopic in $\proj$.
In particular fundamental groups of their complements are isomorphic, 
see Theorem \ref{thepresentation} for a presentation of this group.
In~\cite{holzapfel}, 
the configuration $A_3$ was named the \textit{Apollonius configuration} 
and studied from the orbifold point of view. 
\begin{lemma}\label{qqq}
Let $Q\subset \proj$ be a smooth quadric. 
Then there is a uniformization 
$\psi\,:\,Q\times Q\ok(\proj,2Q)$.
Let $p\in Q$ and put $T^v_p:=\{p\}\times Q$,  
$T^h_p:=Q\times \{p\}$. Then $T_p:=\psi(T^h_p)=\psi(T^v_p)\subset \proj$
 is a line tangent to $Q$ at the point $p\in Q$.
\end{lemma}
\begin{proof}
Since any two smooth quadrics are 
projectively equivalent, it suffices to prove this for a given quadric. 
Consider the $\Z/(2)$-action  defined by  $(x,y)\in\pp\ok (y,x)\in\pp$.
The diagonal $Q=\{(x,x)\,:\,x\in\pline\}$ is fixed under this action.
Let $x=[a:b]\in\pline$ and $y=[c:d]$, then the symmetric polynomials
$\sigma_1([a:b],[c:d]):=ad+bc$, $\sigma_2([a:b],[c,d]):=bd$,
$\sigma_3([a:b],[c:d]):=ac$ are invariant under this action, and the Vi{\'e}te map
$$
\psi:(x,y)\in \pp\longrightarrow [\sigma_1(x,y):\sigma_2(x,y):\sigma_3(x,y)]\in \proj\\
$$
is a branched covering map of degree 2. The branching locus $\subset\proj$
can be found as the image of $Q$. Note that the restriction of 
$\psi$ to the diagonal $Q$ is one-to-one, 
so that one can denote $\psi(Q)$ by the letter
$Q$ again. One has $\psi(Q)=[2ab:b^2:a^2]$  ($[a:b]\in\pline$),
so that $Q$ is a quadric given by the equation $4yz=x^2$. 
One can identify the surface $\pp$ with $Q\times Q$, via the 
projections of the diagonal $Q\subset \pp$.
Let $p\in Q$, and put $T^h_p:=Q\times\{p\}$, $T^v_p:\{p\}\times Q$.
Then $T_p:=\psi(T^h_p)=\psi(T^v_p)\subset\proj$ is a line  tangent to $Q$. 
Indeed, if $p=[a:b]$, then $\psi(T^h_p)$ is parametrized as $[cb+da:db:ca]$ ($[c:d]\in \pline$),
and can be given by the equation $b^2z+a^2y-abx=0$, which shows that  
$T_p$ is tangent to $Q$ at the point $[2ab:b^2:a^2]$.

\end{proof}
Consider the pair $(\proj, aQ+b_1T_1+\dots +b_nT_n)$. 
By Lemma~\ref{tangency} this is an orbifold provided $1/a+1/b_i\geq1/2$. 
Denote
$$
\ca(a;b_1,\dots,b_n):=(\proj,aQ+b_1T_1\dots +b_nT_n)
$$ 
\begin{theorem} \label{firstlifting}
Suppose that $n>1$ and if $n=2$ then $b_1=b_2$. 
Then there is a finite uniformization $\xi:R\times R\ok \ca(2;b_1,\dots,b_n)$,
the Riemann surface $R$ being a uniformization of
$(Q,b_1p_1+\dots+b_np_n)$, where $p_i:=T_i\cap Q$. Moreover;\\
(i) If $n=2$ and $b:=b_1=b_2<\infty$ or $n=3$ and $b_1\1+b_2\1+b_3\1>1$ 
then $R\simeq \pline$. Furthermore, one has $|\orbfg(\ca(2;b,b))|=2b^2$ and
$|\orbfg(\ca(2;b_1,b_2,b_3))|=8[b_1\1+b_2\1+b_3\1-1]\2$ (see Figure 2)\\
(ii) If $n=2$, $b_1=b_2=\infty$ or  $n=3$ and $b_1\1+b_2\1+b_3\1=1$ or $n=4$ and 
$b_1=b_2=b_3=b_4=2$
then $R$ is elliptic. In this case the universal uniformization of $\ca$ is $\C\times\C$ 
(see Figure 1)\\
(iii) Otherwise $R$ is of genus$>1$
and the universal uniformization of $\ca$ is $\pball$. 
\end{theorem}
\begin{proof}
By Lemma~\ref{qqq}, the lifting of  $\ca(2a;b_1,\dots,b_n)$ 
to the uniformization of $\ca(2)=(\proj,2Q)$ is the orbifold (see Figure 1)
$$
\cb(a,b_1,\dots,b_n):=\bigl(Q\times Q, aQ+\sum_{1\leq i\leq n} b_i(T^v_{p_i}+T^h_{p_i})\bigr)
$$
Consequently, there is a covering 
$\psi:\cb(a;b_1,\dots,b_n)\ok \ca(2a;b_1,\dots,b_n)$.
\begin{figure}\label{covering1}
\begin{center}
\resizebox{0.6\hsize}{!}{\includegraphics*{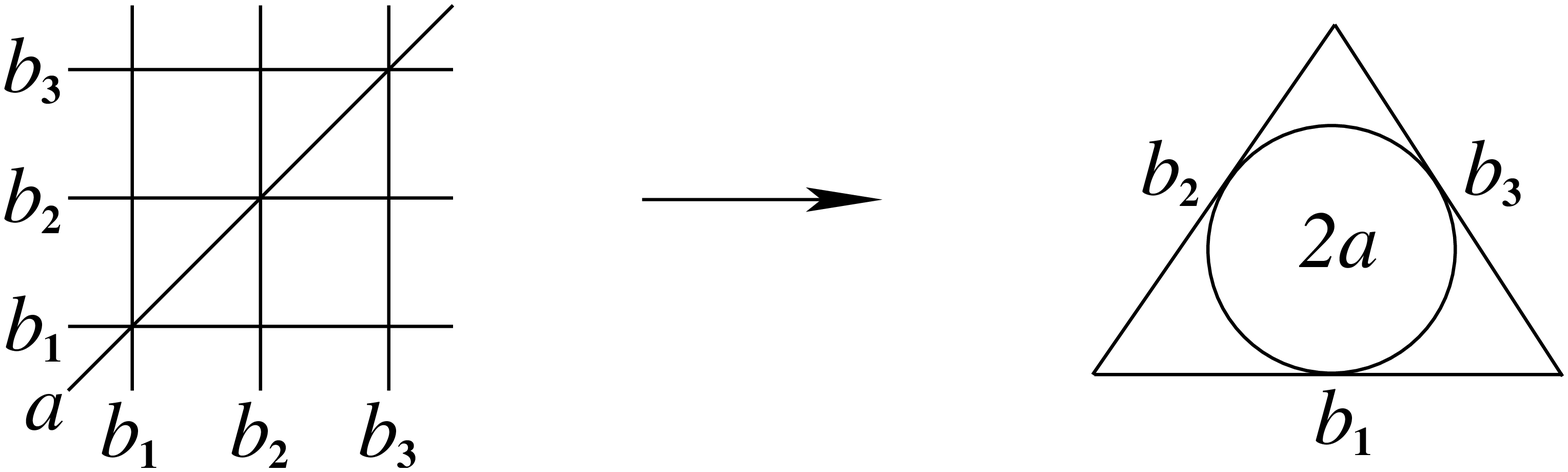}} 
\end{center}
\caption{The covering $\cb(a;b_1,\dots,b_n)\ok \ca(2a;b_1,\dots,b_n)$} 
\end{figure}
Consider the case $a=2$, and denote 
$\cb^\prime(b_1,\dots,b_n):=\cb(1;b_1,\dots,b_n)$.
Then $\cb^\prime=\cs\times\cs$,
where $\cs$ is the one-dimensional orbifold
$\cs(b_1,\dots b_n):=(Q,b_1p_1+\dots +b_np_n)$.
Recall that $Q\simeq\pline$.
Assume that $n>1$ and if $n=2$, then $b_1=b_2$. 
Then $\cs$ admits a finite uniformization $f:R\ok \cs$ by
the Bundgaard-Nielsen-Fox theorem (\cite{bundgaard},\cite{fox2}). 
The Riemann surface $R$ is of genus $0$, $1$ or $>1$ 
according to the conditions stated in the theorem.
It is well known that  in case 
$R\simeq \pline$ the degree of $f$ is the order of the triangle group $\orbfg(\cs)$, 
which is $b^2$ if $n=2$, $b:=b_1=b_2<\infty$ 
and is $2[b_1\1+b_2\1+b_3\1-1]\1$ if $n=3$, $b_1\1+b_2\1+b_3\1>1$. 
\par
The map
$ \zeta:(x,y)\in R\times R\ok (f(x),f(y))\in \cs\times\cs\simeq \cb^\prime$   
is a uniformization of $\cb^\prime$ with 
$\deg(\zeta)=\deg(f)^2$. Since $\zeta$ is compatible 
with the involution on $\cb^\prime$, the composition
$\xi:=\psi\comp\zeta:R\times R\ok \ca$
is a Galois uniformization of $\ca$ with
$\deg(\xi)=\deg(\psi)\deg(\zeta)=2\deg(f)^2.$

\end{proof}
\begin{figure} \label{parabolic}
\begin{center} 
\resizebox{0.15\hsize}{!} {\includegraphics*{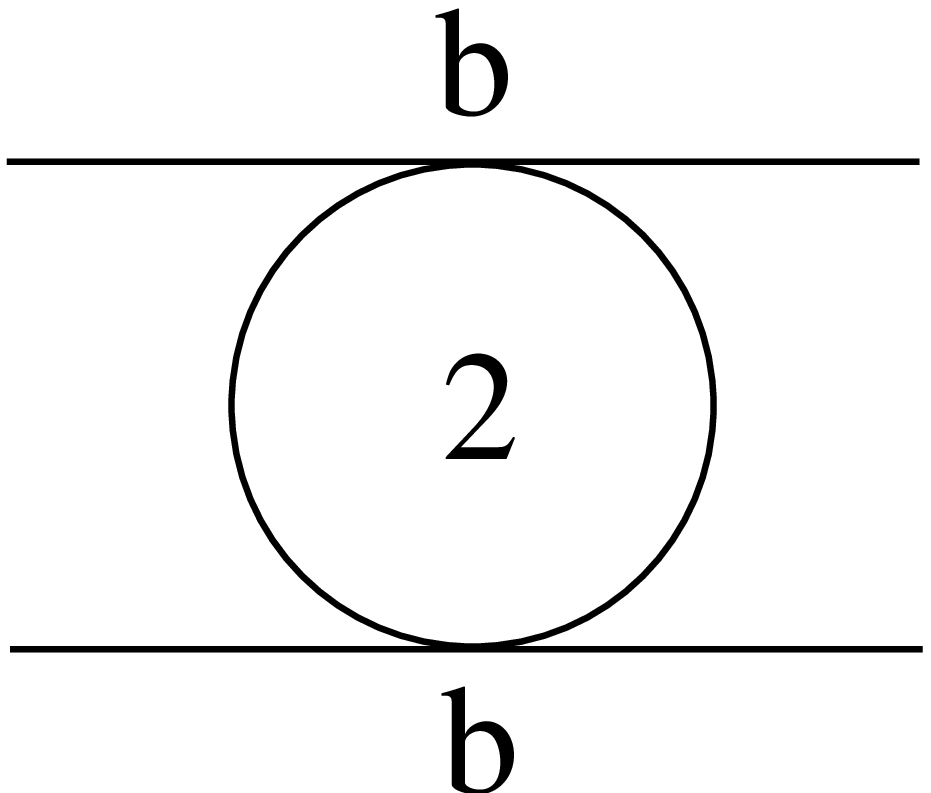}}  
\hspace{3mm} 
\resizebox{0.15\hsize}{!}{\includegraphics*{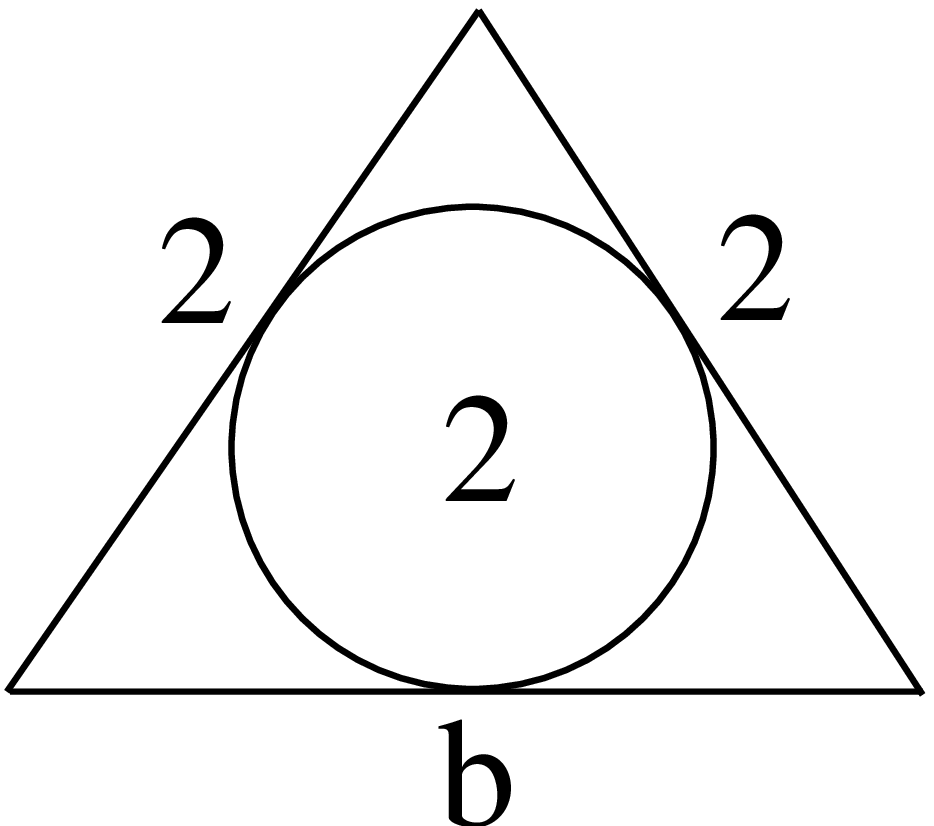}}   
\hspace{3mm} 
\resizebox{0.15\hsize}{!}{\includegraphics*{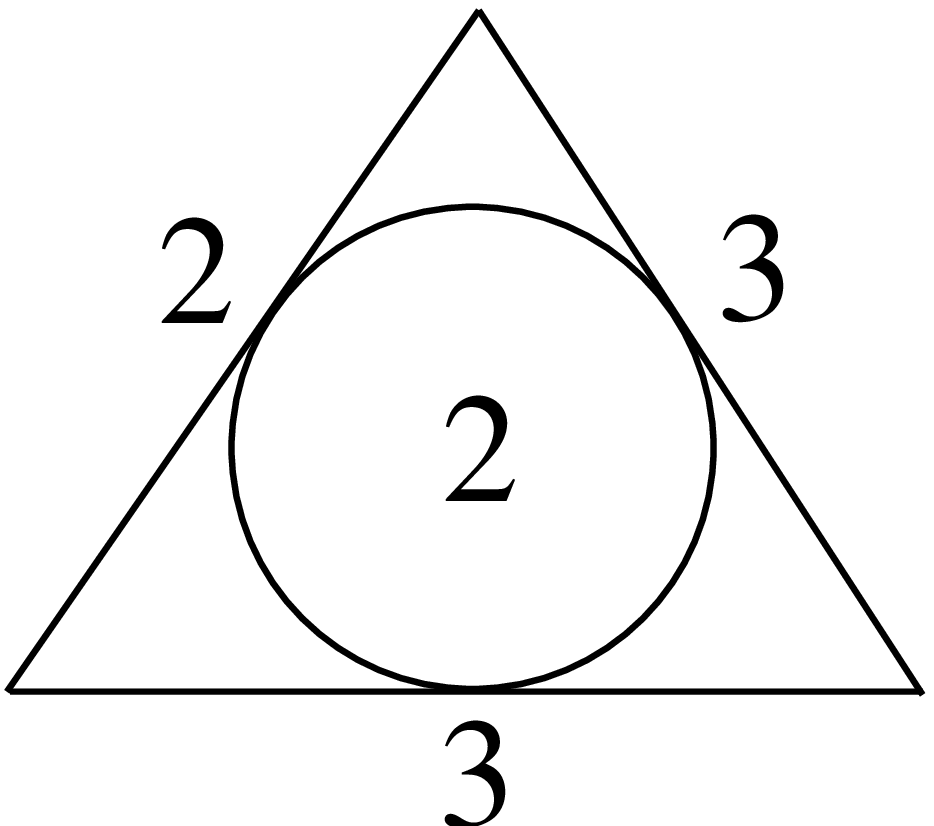}}  
\hspace{3mm} 
\resizebox{0.15\hsize}{!}{\includegraphics*{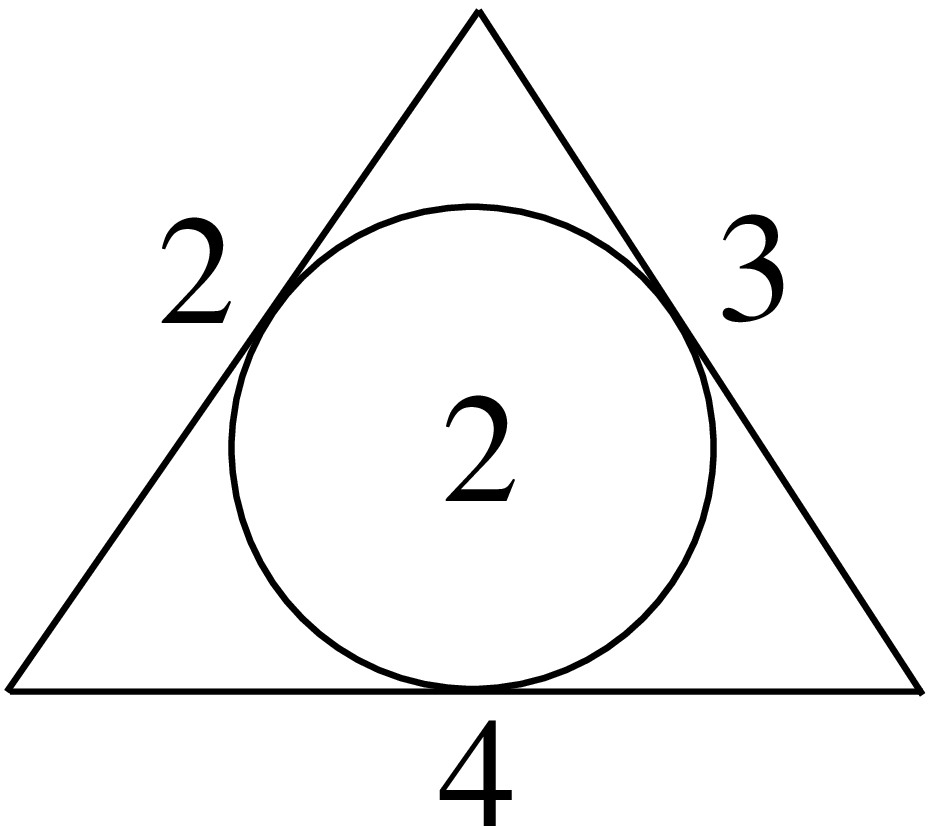}}  
\hspace{3mm} 
\resizebox{0.15\hsize}{!}{\includegraphics*{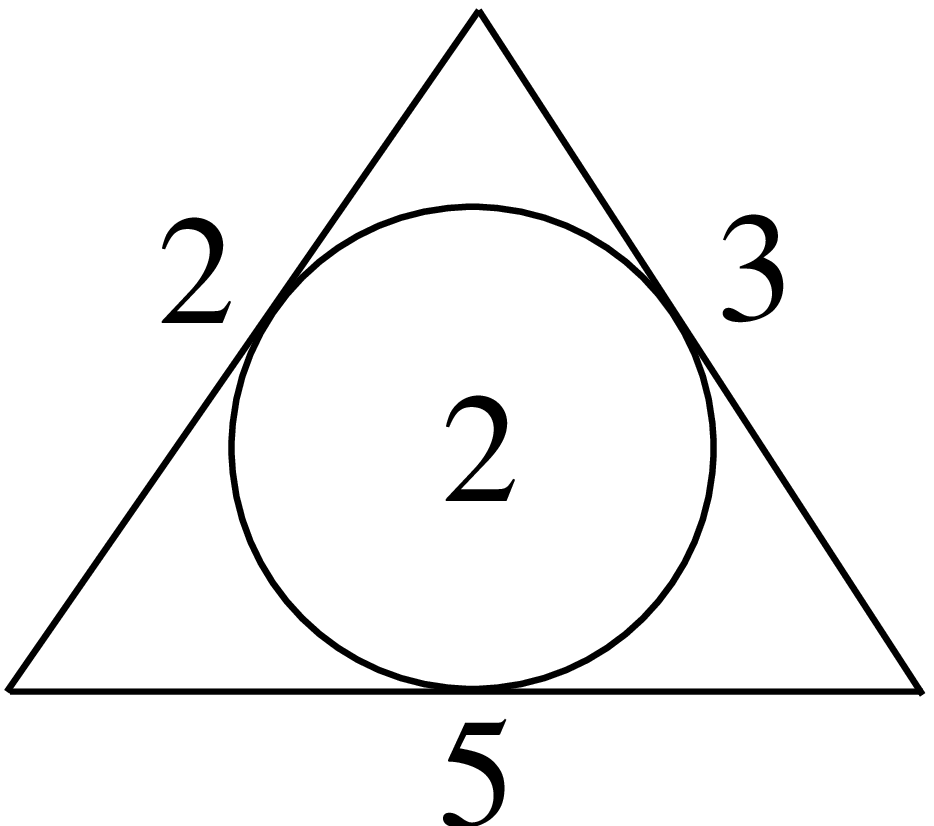}}  
\end{center} 
\caption{Orbifolds $\ca$ uniformized by $\pp$}
\end{figure} 
\begin{figure}\label{elliptic} 
\begin{center}
\resizebox{0.15\hsize}{!}{\includegraphics*{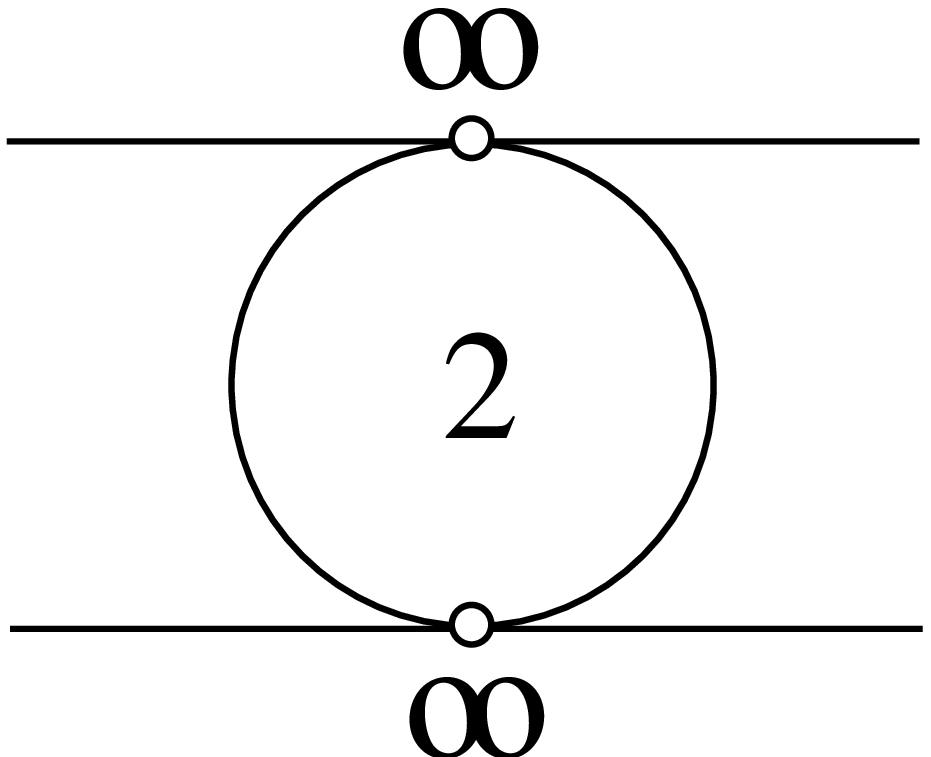}}  
\resizebox{0.15\hsize}{!}{\includegraphics*{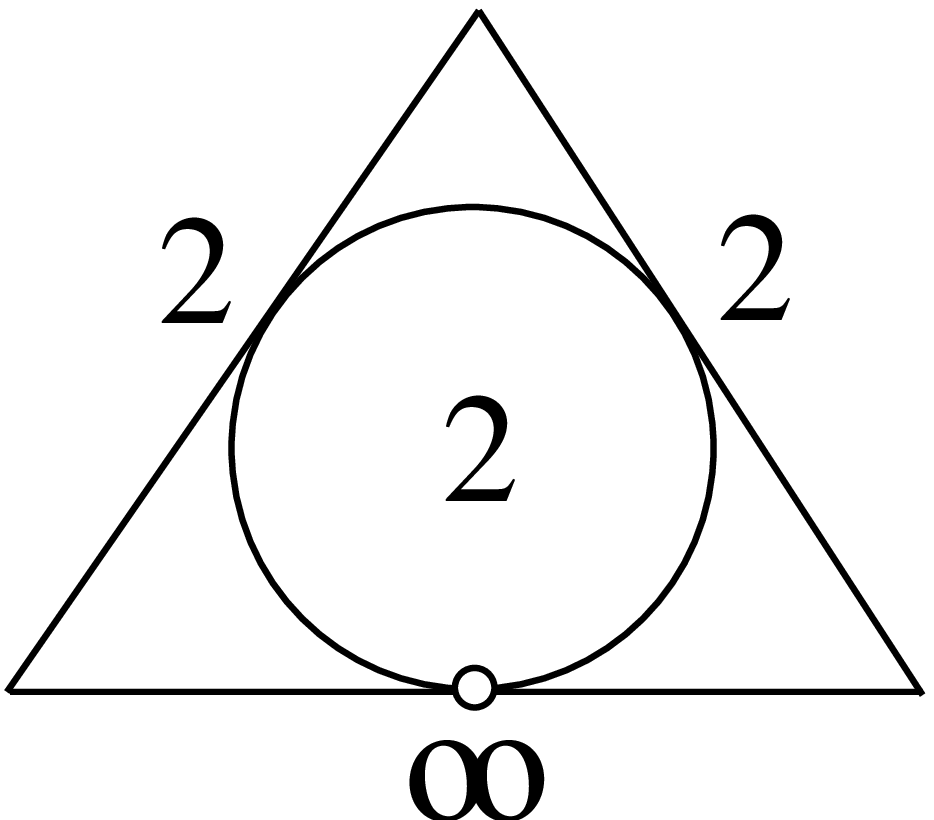}}  
\resizebox{0.15\hsize}{!}{\includegraphics*{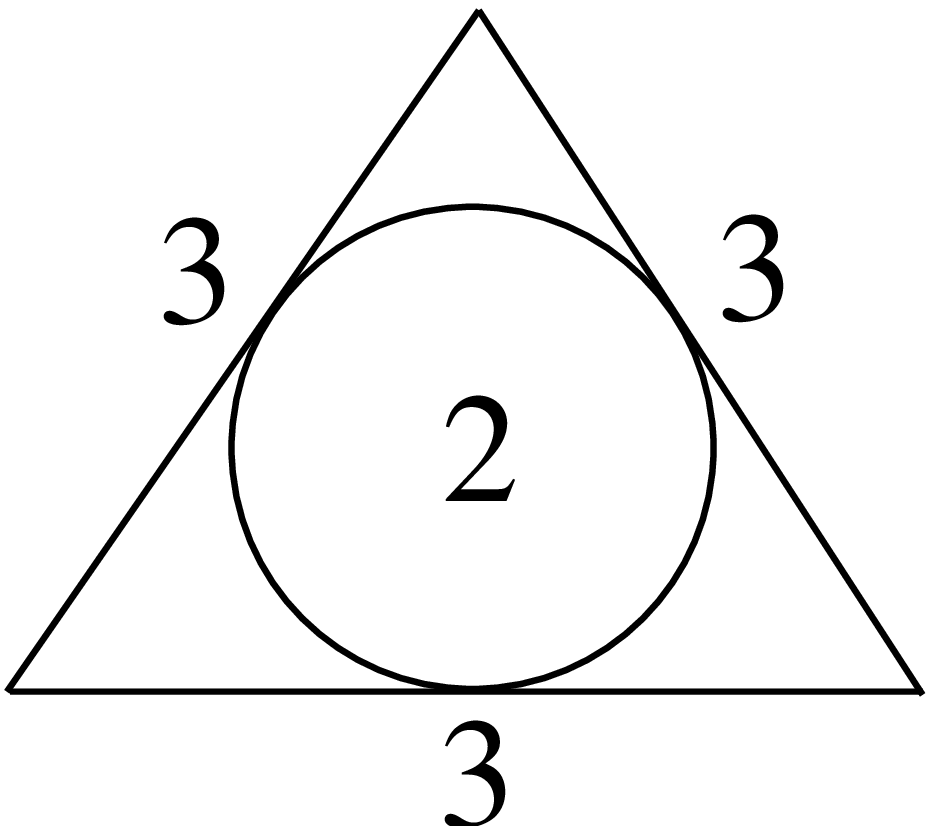}} 
\resizebox{0.15\hsize}{!}{\includegraphics*{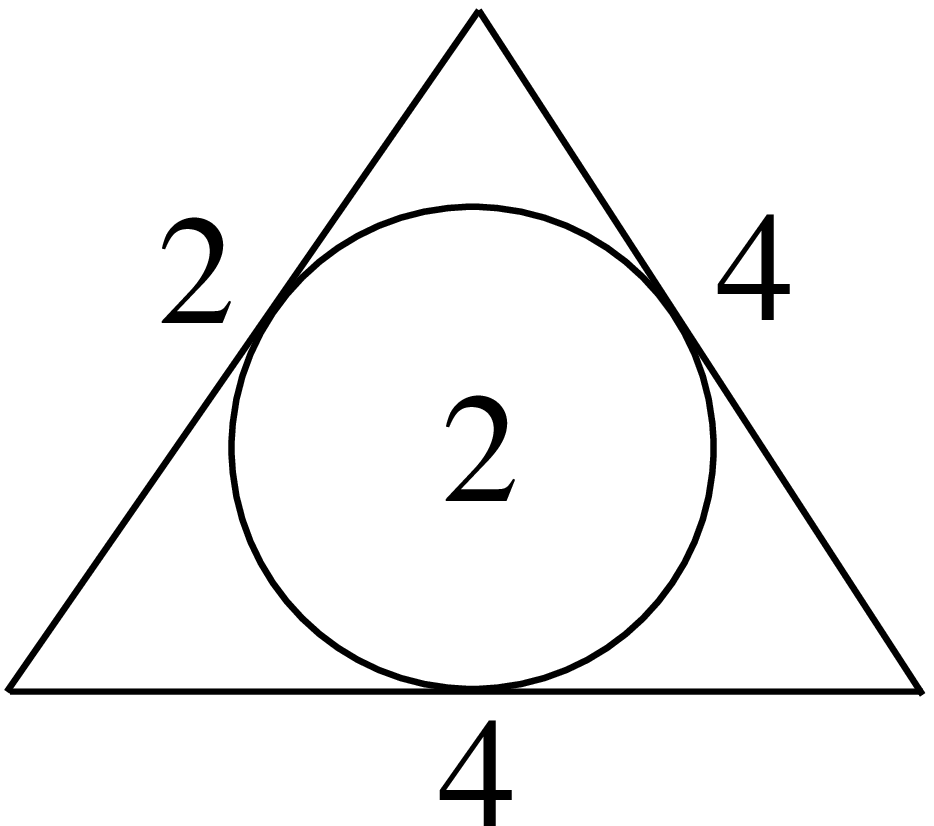}}  
\resizebox{0.15\hsize}{!}{\includegraphics*{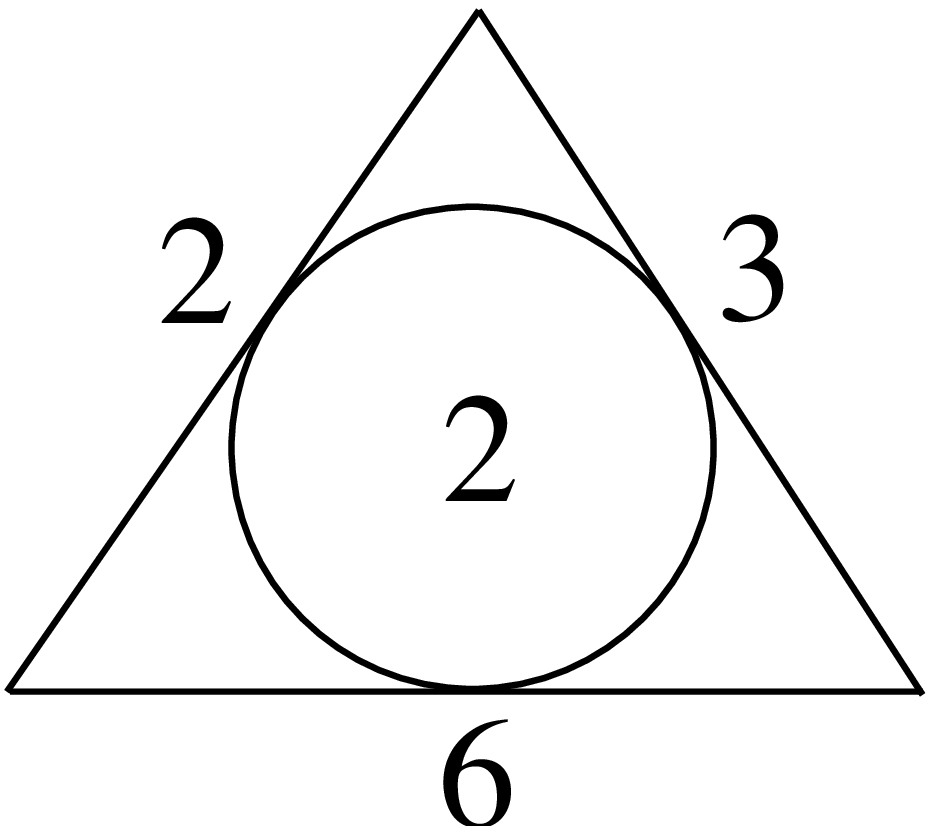}}  
\resizebox{0.13\hsize}{!}{\includegraphics*{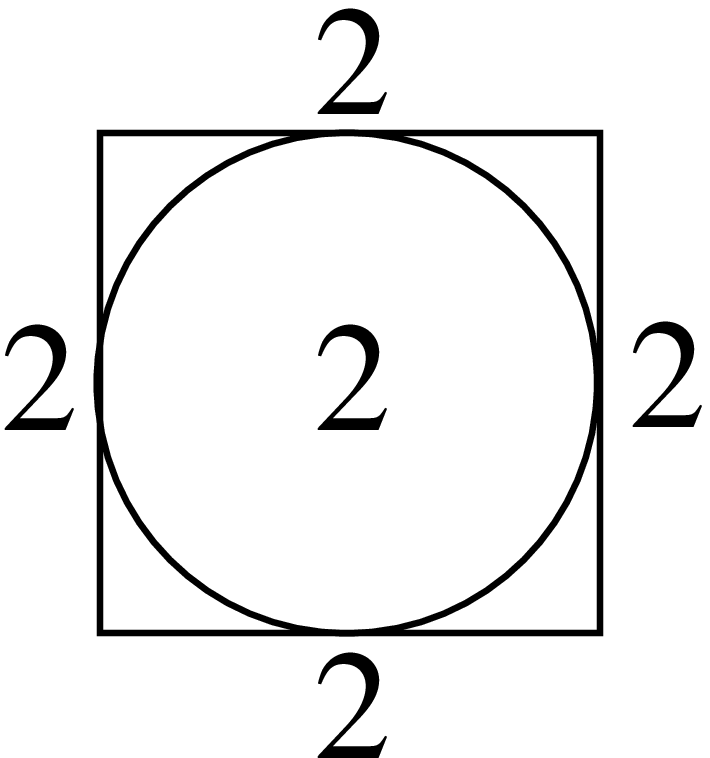}}  
\end{center} 
\caption{Orbifolds $\ca$ uniformized by $\C\times\C$} 
\end{figure}

The orbifolds in Figure~\ref{elliptic} were discovered in~\cite{kaneko2}, where a complete 
classification of the orbifolds $(\proj,D)$ uniformized by $\C\times\C$ was 
also given. Besides these ones  there are two other orbifolds,
which we shall rediscover  in Section~\ref{sectioncoverings}. 

\subsection{Chern numbers of the orbifolds $\ca(a;b_1,\dots,b_n)$}\label{sectioncherns} 
By definition, the Chern numbers of $\ca(a;b_1,\dots,b_n)$ are given by 
$$
c_1^2(\ca)=\bigl[n-1-{2}{a}\1-\sum_{1\leq i\leq n} {b_i}\1\bigr]^2
$$
\begin{eqnarray*}
e(\ca)=\frac{(n-1)(n-2)}{2}+{(2-n)}{a\1}+\sum_{1\leq i\leq n}{(2-n)}{b_i\1} \qquad\qquad\qquad \qquad 
\\
\qquad \qquad \qquad  +\sum_{1\leq i\neq j\leq n}{b_i\1b_j\1}+\frac{1}{2}
\sum_{1\leq i\leq n}\left[{b_i}\1+{a}\1-{2}\1\right]^2
\end{eqnarray*}
\begin{proposition}\label{parabolics}
For the orbifold $\ca(a;b_1,\dots,b_n)$, one has \\
(i) $2e=c^2$ if and only if  $a=2$ or 
$n=4$, $b_1=b_2=b_3=b_4=2$ or 
$n=3$, $(a;b_1,b_2,b_3)=(3;2,3,4)$.\\
(ii) $e=c^2=0$ if and only if $n=2$, $a=2$; $b_1=b_2=\infty$ or $n=3$, $a=2$; 
$b_1\1+b_2\1+b_3\1=1$ or $n=4$, $a=2$; $b_1=b_2=b_3=b_4=2$. \\
(iii) $3e=c_1^2>0$ if and only if $n=3$ and
$(a;b_1,b_2,b_3)$ is one of $(4;4,4,4)$, $(3;3,4,4)$, $(3;6,6,2)$ or $(3;6,3,3)$.\\
(iv) $c_1^2=0$ and $e>0$ if and only if $n=2$, {and $(a;b_1,b_2)$ is one of $(4;4,4)$, 
$(3;6,6)$, $(6;3,3)$} or $n=3$ and 
$(a;b_1,b_2,b_3)$ is one of $(4;2,2,2)$ or $(3;3,2,2)$. 
\end{proposition}
\begin{proof}
Put $\beta=\sum_{1\leq i\leq n}b_i\1$  and $\alpha:=\beta-a\1$.
We shall prove the following splitting formulas for the difference of the Chern numbers of $\ca$,
the claims of the proposition follows easily from these formulas. 
\begin{equation}\label{firstprop}
2(2e-c_1^2)(\ca)=( {a}\1-2\1 )
\bigl[n(2a\1+{3})-4\beta-4(2a\1+{1})\bigr]
\end{equation}
\begin{equation}\label{secondprop}
8(3e-c_1^2)(\ca)=(2\alpha+5-2n)^2+3(n-3)\left(2{a}\1-1\right)^2
\end{equation}
To prove (\ref{firstprop}), first note that 
$$
c_1^2(\ca)=(n-1)^2+4a\2+\beta^2-4(n-1)a\1-2(n-1)\beta+4a\1\beta 
$$
and
\begin{eqnarray*}
2e(\ca)=(n-1)(n-2)+2(2-n)a\1+2\beta(2-n) +2\sum_{1\leq i\neq j\leq n}b_i\1b_j\1+\\
\sum_{1\leq i\leq n}\left[b_i\2+a\2+2\2+2b_i\1a\1-b_i\1-a\1\right]^2
\end{eqnarray*}
$$
=(n-1)(n-2)+\frac{n}{4} +2(2-n)a\1+2\beta(2-n)+\beta^2+na\2-na\1+2\beta a\1-\beta 
$$
which gives
\begin{equation}\label{thelastexp}
(2e-c_1^2)(\ca)=1-\frac{3n}{4}+na\1+\beta+(n-4)a\2-2a\1\beta
\end{equation}
The formula~(\ref{firstprop}) is easily seen to be equivalent to~(\ref{thelastexp}).
The formula~(\ref{secondprop}) is proved similarly, by using the above expressions for 
$e(\ca)$ and $c_1^2(\ca)$.
\end{proof}
The orbifolds (i)-(ii) in Proposition~\ref{parabolics} 
were shown to be uniformizable in Theorem~\ref{firstlifting}.
Orbifolds $\ca(a;2,2,2,2)$ will be shown to be uniformizable in
Proposition~\ref{remainings} below.
We don't know if the orbifold $\ca(3;2,3,4)$ in 
$(i)$ is uniformizable. The uniformizability by $\ball$ of the orbifolds (iii) 
follows from the KNS theorem~\cite{kobayashi} or from~\cite{holzapfel3}, 
where it is also shown that the corresponding lattices are arithmetic. 
{As for the orbifolds (iv), it will be shown elsewhere that the orbifolds 
$\ca(4;4,4)$, $\ca(3;6,6)$ and $\ca(6;3,3)$ are not uniformizable, and that 
the orbifold $\ca(3;3,2,2)$ is uniformized by a K3 surface.} 
The universal uniformization of the orbifold $\ca(4;2,2,2)$ is a K3 surface,
and the group $\orbfg(\ca(4;2,2,2))$ is finite of order 256, see Proposition~\ref{fgs}.

\begin{figure}\label{hyperbolic} 
\begin{center}
\resizebox{0.15\hsize}{!}{\includegraphics*{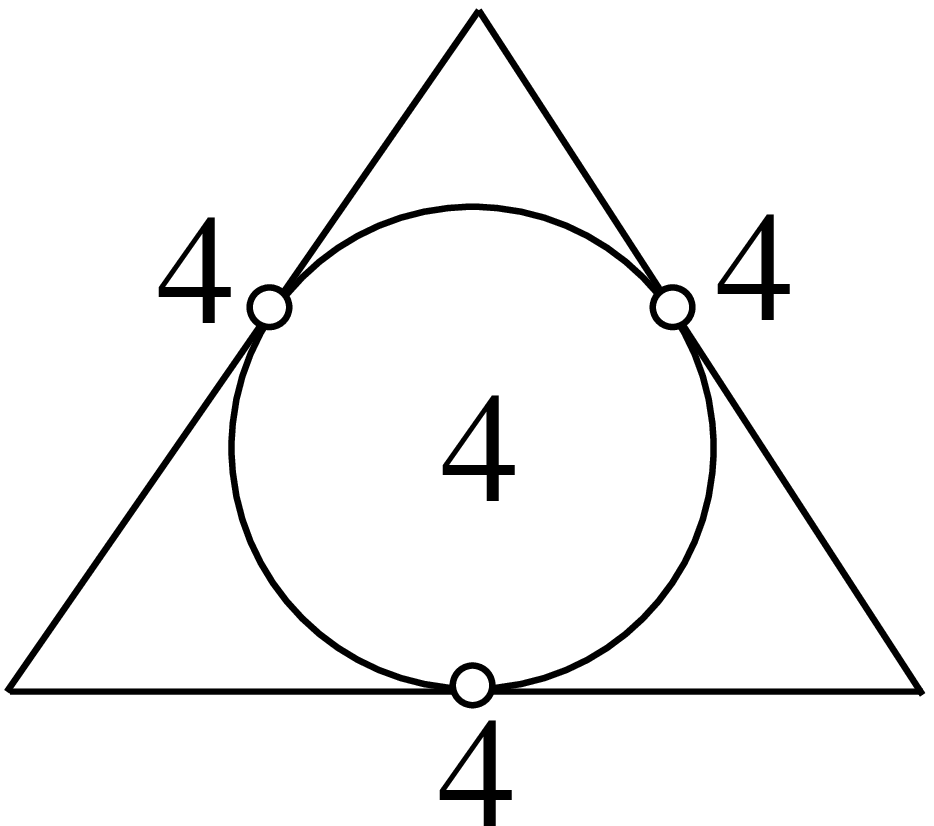}}  
\hspace{4mm} 
\resizebox{0.15\hsize}{!}{\includegraphics*{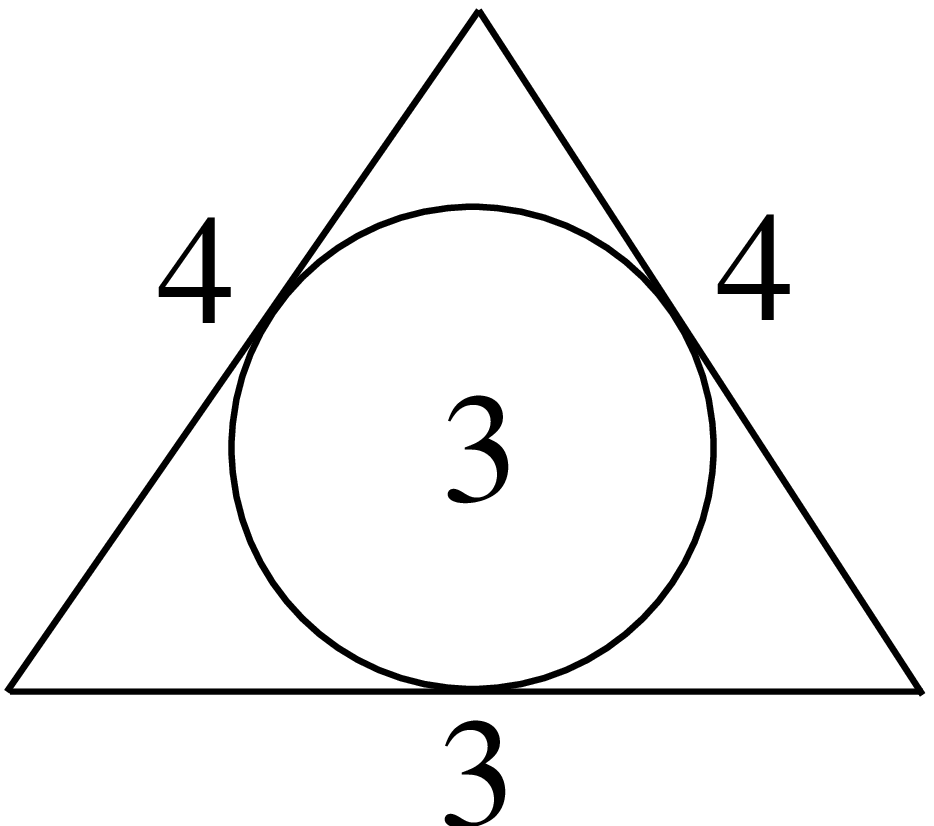}}  
\hspace{4mm} 
\resizebox{0.15\hsize}{!}{\includegraphics*{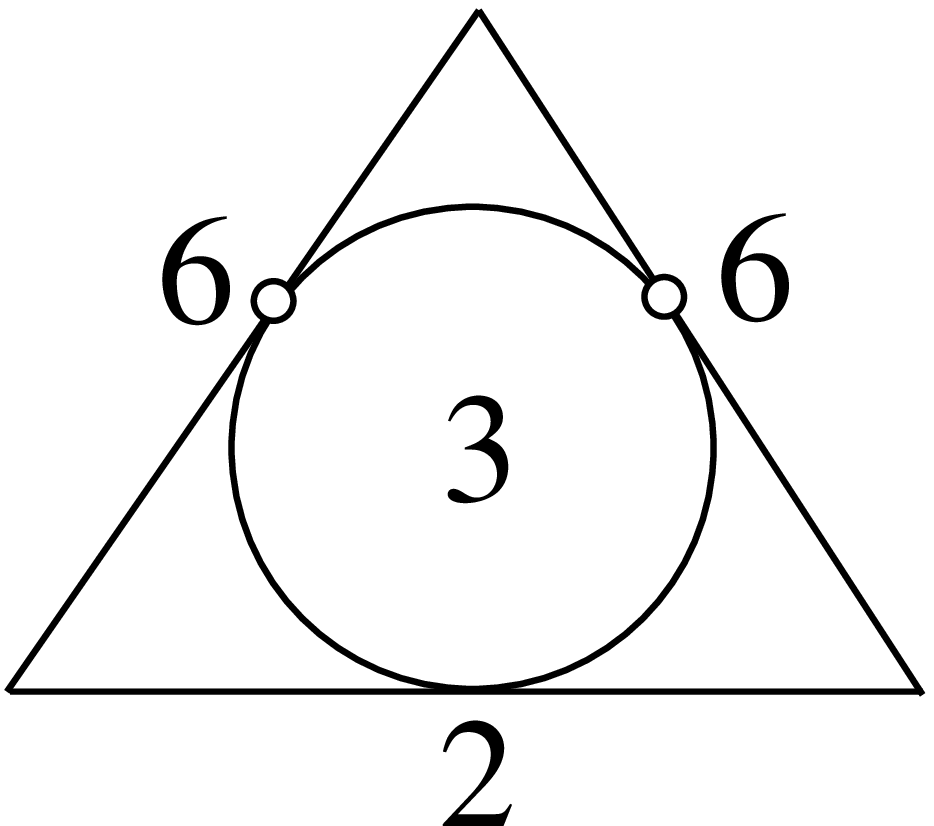}} 
\hspace{4mm} 
\resizebox{0.15\hsize}{!}{\includegraphics*{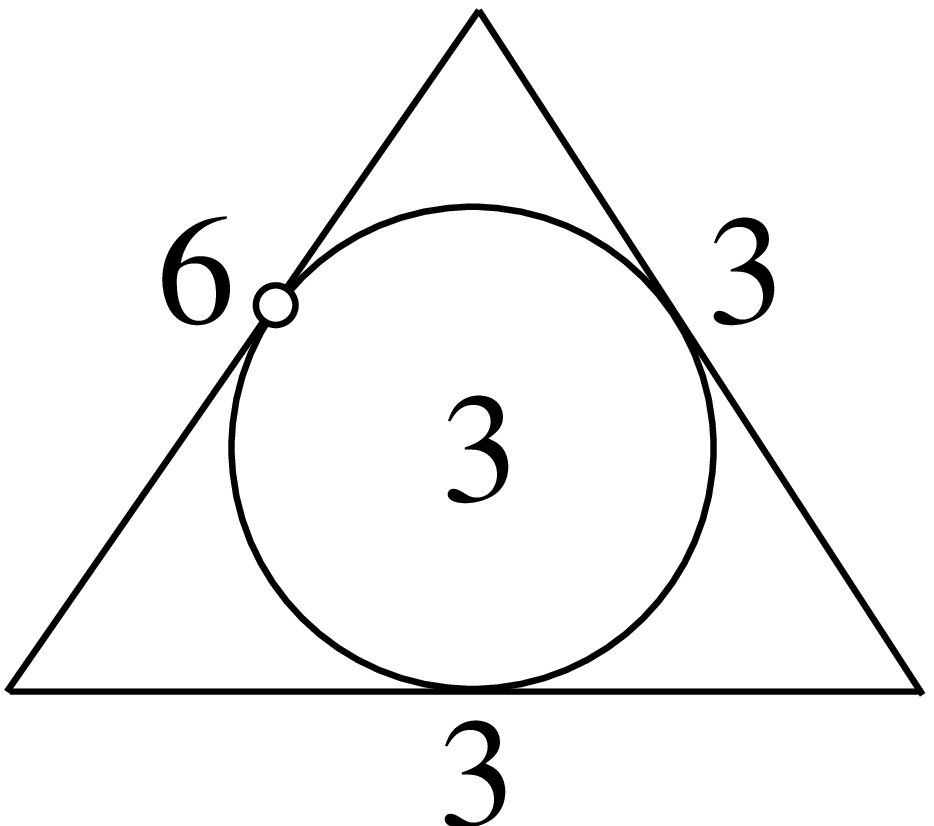}}  
\end{center} 
\caption{Orbifolds $\ca$ uniformized by $\ball$} 
\end{figure} 
\begin{figure}\label{calabiyau} 
\begin{center}
\resizebox{0.15\hsize}{!}{\includegraphics*{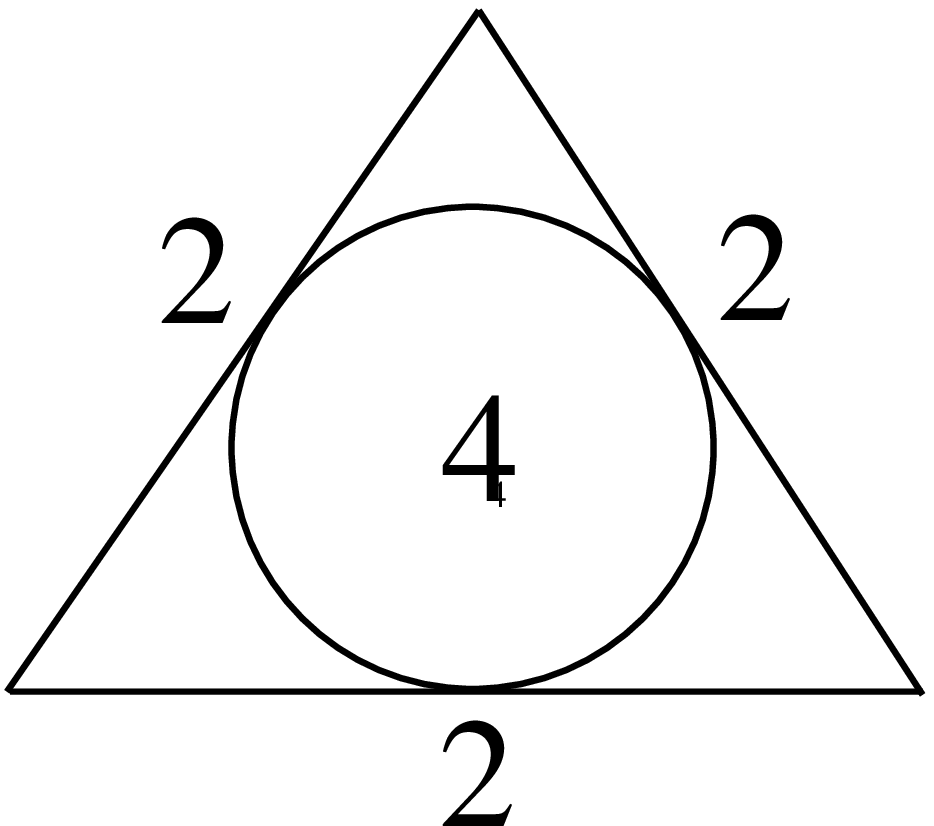}}  
\hspace{6mm} 
\resizebox{0.15\hsize}{!}{\includegraphics*{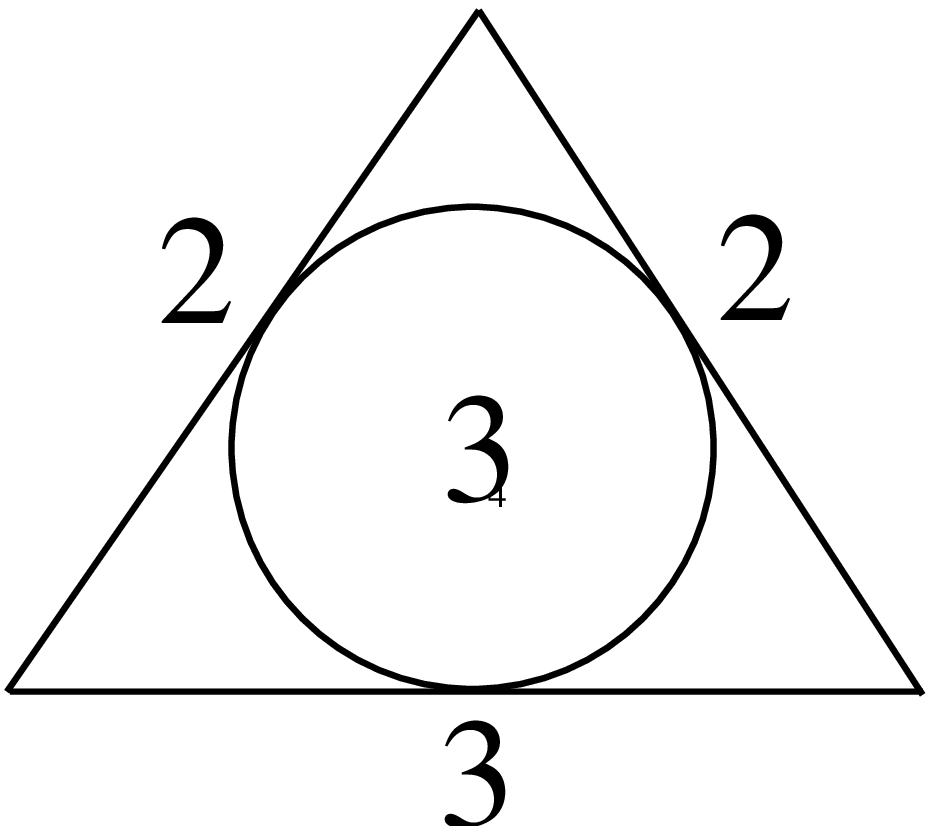}} 
\end{center} 
\caption{Orbifolds $\ca$ uniformized by  K3 surfaces} 
\end{figure} 

\subsection{Another orbifold covering of $\ca$}
Let $\ca=(\proj, aQ+\sum_{1\leq i\leq n}b_iT_i)$ 
be an orbifold over the Apollonius configuration with $n\geq 3$.
Suppose without loss of generality that 
the lines $T_1$, $T_2$ and $T_3$ are defined by the equations 
$x=0$, $y=0$, and $z=0$.
Let $k\in \N$ and consider the orbifold $(\proj, kT_1+kT_2+kT_3)$. 
This orbifold admits the uniformization
$\phi_k:[x,y,z]\in\proj\ok [x^k,y^k,z^k]\in\proj$.
Denote by $X$, $Y$ and $Z$ the lines $\phi_k\1(T_1)$, $\phi_k\1(T_2)$ 
and $\phi_k\1(T_3)$ respectively.
For $i\geq 4$,  $T_i$ is given by the equation $\alpha x+\beta y+\gamma z=0$, then 
$T_{i,k}:=\phi_k\1(T_i)$ 
is the Fermat curve given by the equation $\alpha x^k+\beta y^k+\gamma z^k=0$.
\par
The equation of the quadric $Q$, tangent to the lines $T_1$, $T_2$ and $T_3$ 
can be written in the form
\begin{equation}\label{quadriceq}
Q: \alpha\pr x^{1/2}+\beta \pr y^{1/2}+\gamma \pr z^{1/2}=0,
\end{equation}
where $\alpha\pr\beta\pr\gamma\pr\neq 0$.
This shows that for $k$ odd, $Q_k:=\phi_k\1(Q)$ is an irreducible curve given by the equation
$ \alpha\pr x^{k/2}+\beta \pr y^{k/2}+\gamma \pr z^{k/2}=0$, and for $k=2m$ even, $\phi_k\1(Q)$ consists of four 
Fermat curves $Q_{1,k}$, $Q_{2,k}$, $Q_{3,k}$, $Q_{4,k}$ given by the equations 
$\alpha\pr x^m\underline{+} \beta \pr y^m \underline{+}\gamma \pr z^m$. 
\par
For $k$ odd, define the orbifolds
$$
\cc_k(a;e,f,g;b_4,\dots,b_n):=(\proj, aQ_k+eX+fY+gZ+\sum_{i=4}^n b_iT_{i,k}),
$$
and for even $k=2m$, define the orbifolds 
\begin{eqnarray*}
\cc_{k}(a_1,a_2,a_3,a_4;e,f,g;b_4,\dots,b_n):=\qquad\qquad\qquad\qquad\qquad\qquad\qquad\qquad\\
\qquad\qquad\qquad\qquad (\proj, \sum_{j=1}^4 a_jQ_{j,k}+eX+fY+gZ+\sum_{i=4}^n b_iT_{i,k})
\end{eqnarray*}
With these notations, one has the following obvious lemma:
\begin{lemma}\label{secondlifting}
(i) For $k$ odd, there is an orbifold covering
$$
\phi_k:\cc_k(a;e,f,g;b_4,\dots,b_n)\ok\ca(a;ke,kf,kg;b_4,\dots b_n).
$$
(ii) For $k$ even, there is an orbifold covering 
$$
\cc_{k}(a,a,a,a;e,f,g;b_4,\dots,b_n)\ok\ca(a;ke,kf,kg;b_4,\dots b_n).
$$
\end{lemma}
The covering map $\phi_2$ and the orbifold $\cc_2$ are particularly interesting. 
{In this case, $L_1:=Q_{1,2}$, $L_2:=Q_{2,2}$, $L_3:=Q_{3,2}$, $L_4:=Q_{4,2}$ 
are four lines given by the equations
$\alpha\pr x\underline{+} \beta \pr y \underline{+}\gamma \pr z$.}
Since $abc\neq 0$, the points $L_i\cap L_j$ lie on the smooth points of $xyz=0$.
The curves $T_{i,2}$ are smooth quadrics for $i\geq 4$.
Each $T_{i,2}$ has the four lines $L_1$, $L_2$, $L_3$, $L_4$ as tangents.
In particular, for $n=3$ one has
the lines $L_i$ ($i\in\{1,2,3,4\}$) and the lines $X$, $Y$, $Z$,
forming an arrangement of 7 lines with 6 triple points and 3 nodes. 
\begin{figure}\label{covering2}
\begin{center}
\resizebox{0.6\hsize}{!}{\includegraphics*{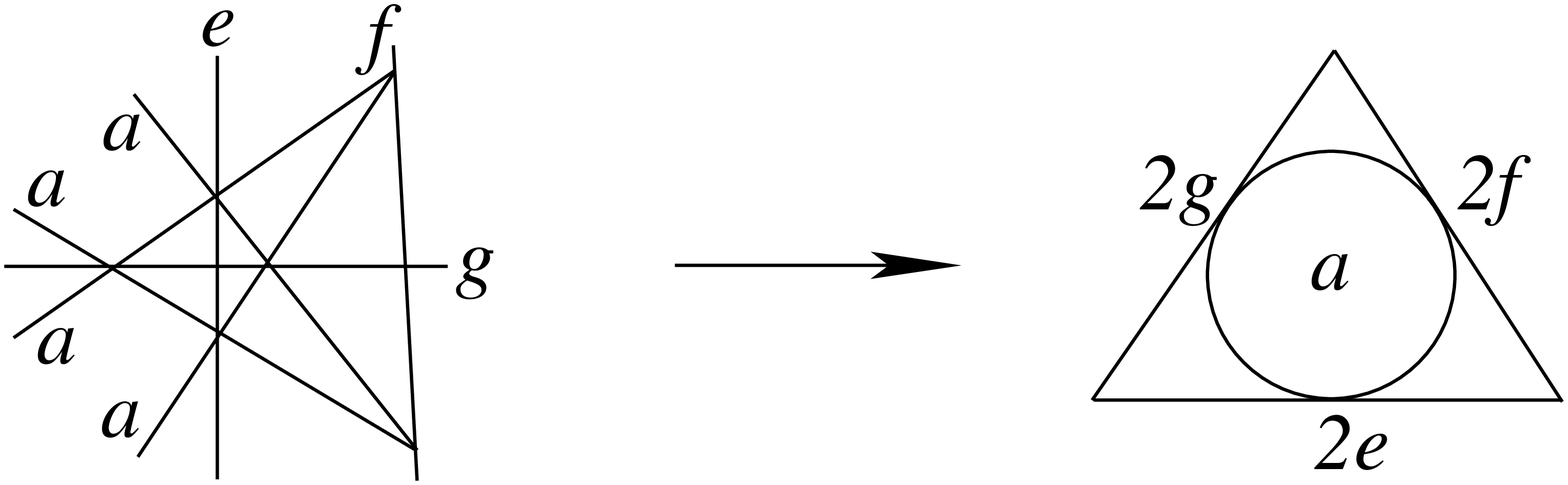}} 
\end{center}
\caption{The covering $\cc_2(a,a,a,a;e,f,g)\ok \ca(a;2e,2f,2g)$} 
\end{figure}
For $n=4$, the curve $T_{4,2}$ is a quadric, and the lines $L_i$ ($i\in\{1,2,3,4\}$)
are tangent to this quadric,  in other words $\cc_2(b,b,b,b;1,1,1;2)\simeq \ca(2;b,b,b,b)$.
This proves the following proposition.
\begin{proposition}\label{remainings}
There is a covering $\phi_2:\ca(2;b,b,b,b)\ok\ca(b;2,2,2,2)$.
Hence the orbifolds $\ca(b;2,2,2,2)$ are uniformized by $\pball$ for $b>2$ .
\end{proposition}
For the orbifolds $\ca(a;2,2,2)$ one has
\begin{proposition}\label{fgs}
The group $\orbfg(\ca(a;2,2,2))$ is finite of order $4a^3$.
There is a universal uniformization $M_a\ok \ca(a;2,2,2)$ 
such that 
$$
e(M_a)=a(a^2-4a+6),\quad c_1^2(M_a)=a(4-a)^2.
$$ 
In particular, $\ca(4;2,2,2)$ is uniformized by a $K3$ surface.
\end{proposition}
\begin{proof}
The orbifold $\cc_2(a,a,a,a;1,1,1)$  admits a universal uniformization of degree $a^3$ since its locus  
consists of the lines $L_i$ $(1\leq i\leq 4)$ in general position and
$\orbfg(\cc_2(a,a,a,a;1,1,1))\simeq \Z/(a)\oplus\Z/(a)\oplus\Z/(a)$.
Composing this with the covering $\cc_2(a,a,a,a;1,1,1)\ok \ca(a;2,2,2)$ gives
 the desired result. 
\end{proof}

\section{Coverings of the orbifolds uniformized by $\pball$.}
\label{sectioncoverings}
\begin{theorem}\label{coverings}
(i) The orbifolds  
$\cc_2(2,2,2,2;e,1,1)$ are uniformized by $\pp$ (see Figure~7)\\ 
(ii) The orbifolds  $\cc_2(2,2,2,2;2,2,1)$ and $\cc_3(2;1,1,1)$ are 
uniformized by $\C\times\C$ (see Figure~7)\\
(iii) Otherwise, $\cc_k(2,2,2,2;e,f,g;b_4,\dots,b_n)$ is  uniformized by $\pball$  
(see Figure~8).
\end{theorem} 
\begin{proof}
By Lemma~\ref{secondlifting}, for $k$ even there is a covering
$$
\cc_k(2,2,2,2;e,f,g;b_4,\dots,b_n)\ok\ca(2;ek,fk,gk;b_4,\dots,b_n),
$$
and for $k$ odd there is a covering
$$
\cc_k(2,2,2,2;e,f,g;b_4,\dots,b_n)\ok\ca(2;ek,fk,gk;b_4,\dots,b_n).
$$
On the other  hand, by Theorem~\ref{firstlifting}, the orbifolds
$\ca(2;ek,fk,gk,b_4,\dots,b_n)$ are uniformized by $\pp$, $\C\times\C$, or 
$\pball$ according to the conditions stated in the theorem. 

\end{proof}

\begin{figure}
\label{unifbyp1p1} 
\begin{center}
\resizebox{0.17\hsize}{!}{\includegraphics*{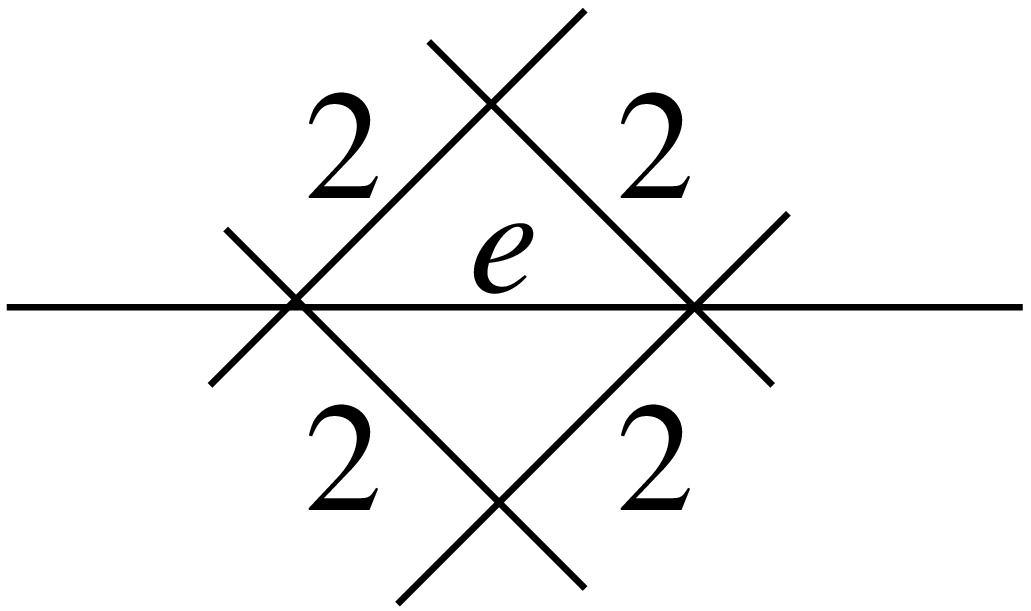}}  
\hspace{1cm}
\resizebox{0.17\hsize}{!}{\includegraphics*{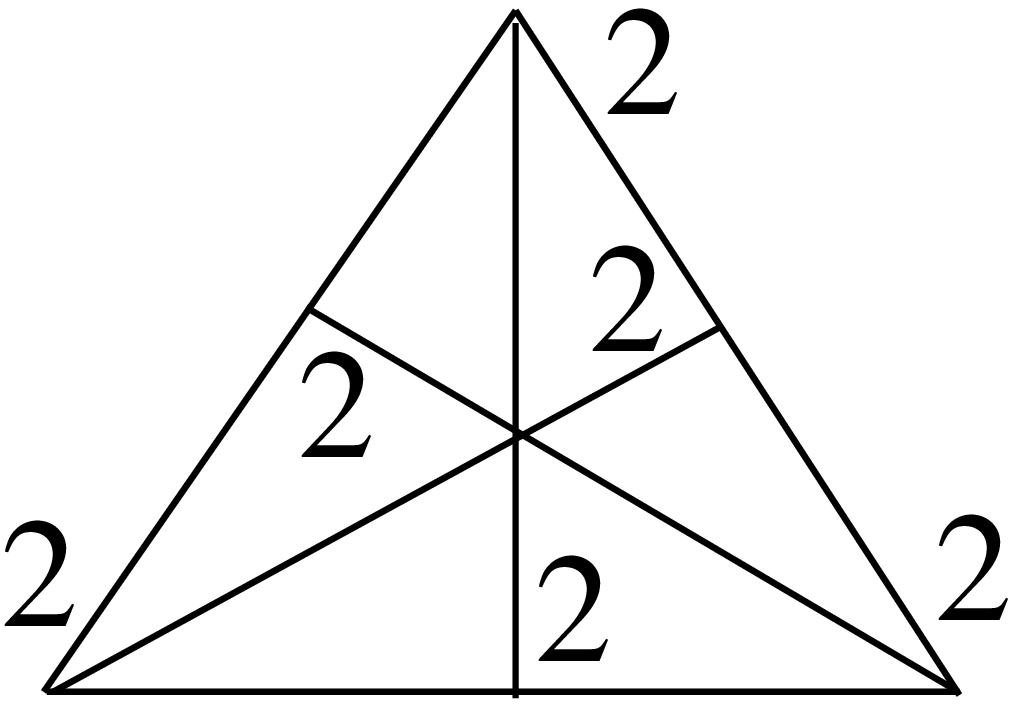}}  
\caption{Orbifolds $\cc_2$ uniformized by $\pline\times\pline$ and by $\C\times\C$} 
\end{center}
\end{figure} 

\begin{figure}
\label{unifbyb1b1}
\begin{center} 
\resizebox{0.17\hsize}{!}{\includegraphics*{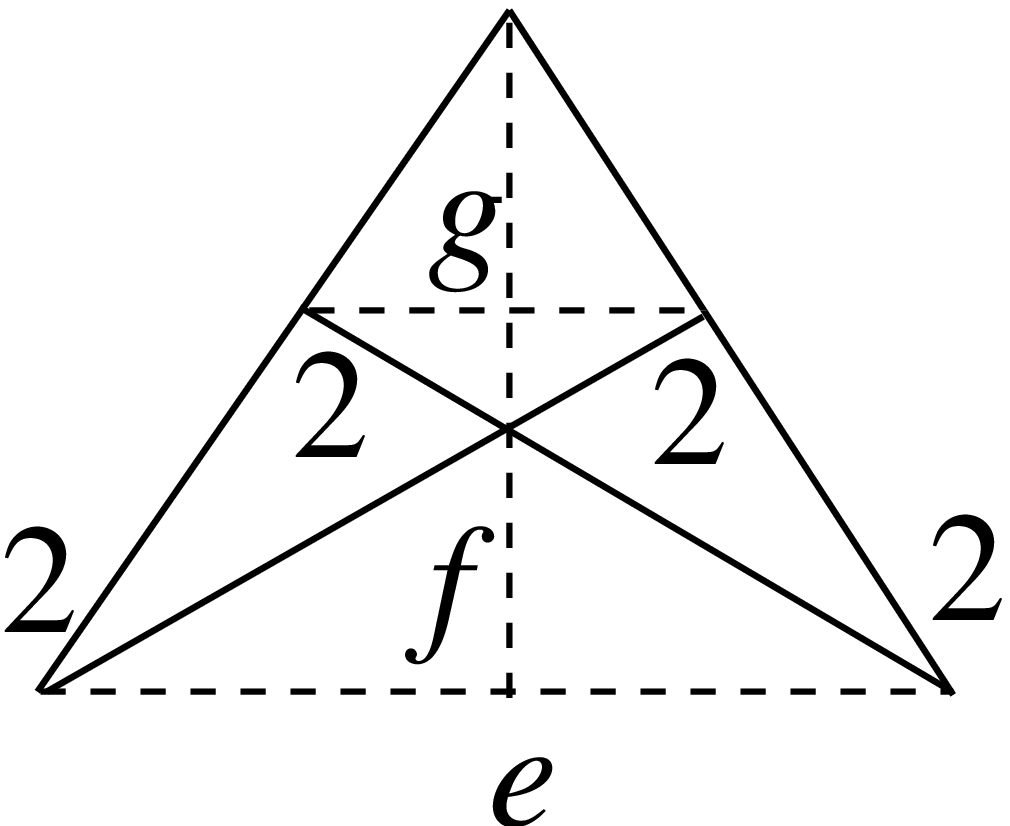}}  
\hspace{1cm}
\resizebox{0.17\hsize}{!}{\includegraphics*{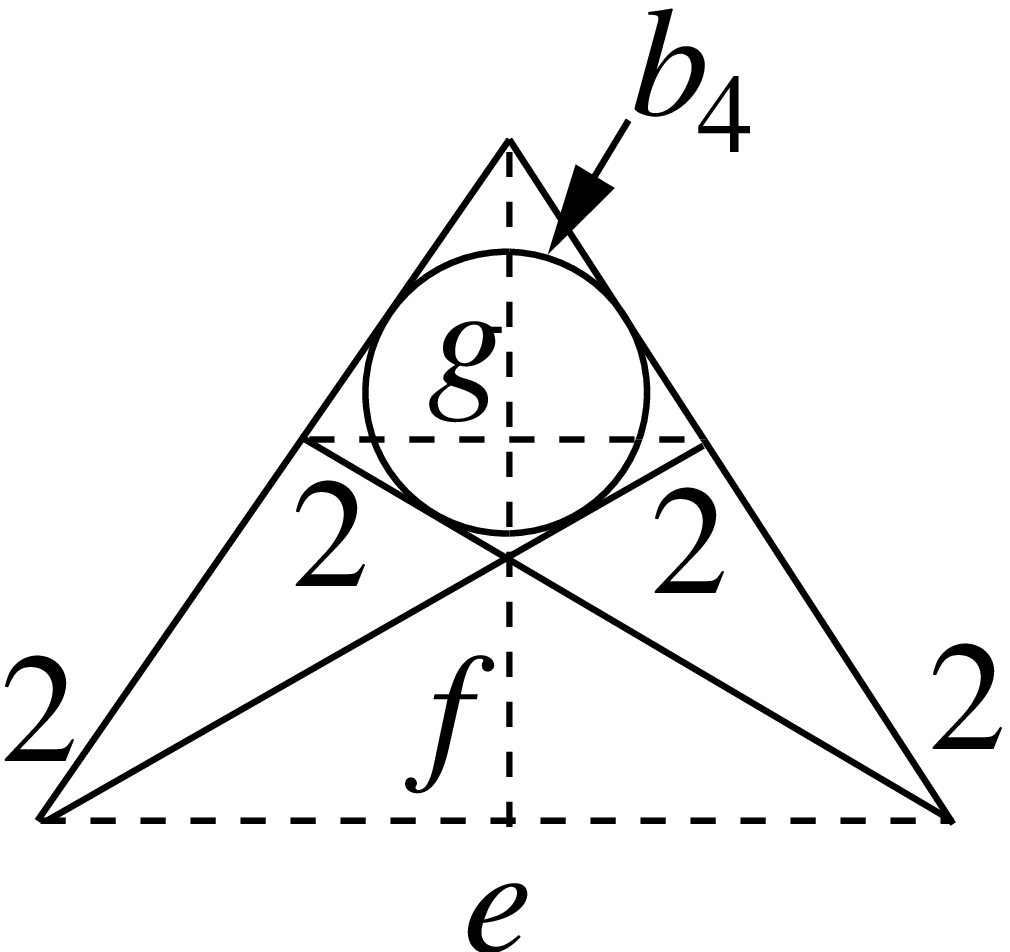}}  
\end{center}
\caption{Some orbifolds $\cc_2$ uniformized by $\pball$} 
\end{figure} 
Consider the line arrangement $\{L_1,L_2,L_3,L_4,X,Y,Z\}$.
The  lines $\{L_1, $ $L_2,L_3,$ $L_4, X, Y\}$ 
forms a complete quadrilateral, i.e. an arrangement of six lines 
with four triple points and three nodes. 
The triple points of the complete quadrilateral 
can be given as $L_1\cap L_2\cap Y$, $L_2\cap L_3\cap X$, $L_3 \cap L_4\cap Y$
and $L_4\cap L_1\cap X$, in this case the nodes are $X\cap Y$, $L_1\cap L_3$ and $L_2\cap L_4$. 
The line $Z$ passes through the nodes $L_1\cap L_3$ and $L_2\cap L_4$.
Take the lines $\{L_1,L_2, X\}$, which do not meet at a triple point, and consider 
the uniformization map $\gamma_2:\proj\ok (\proj, mL_1+mL_2+mX)$.  
Then  $\gamma_2\1(L_3)$ (or $\gamma_2\1(L_4)$ or $\gamma_2\1(Y)$) will consist of $m$ lines,
forming an arrangement of $3m$ lines with $m^2$ triple points and three points
of multiplicity $m$. The curve $Z\pr:=\gamma_2\1(Z)$ will be a smooth Fermat curve of degree $m$.
In the particular case where $m=2$, let
$$
\{L_{3}\pr,L_{3}\prpr\}:=\gamma_2\1(L_3), \quad  \{L_{4}\pr,L_{4}\prpr\}:=\gamma_2\1(L_4),\quad \{Y\pr,Y\prpr\}:=\gamma_2\1(Y),
$$
and define the orbifold $\cd(a_1,a_2;a_3\pr,a_3\prpr;a_4\pr,a_4\prpr; e;f\pr,f\prpr;g)$ as 
$$
(\proj, a_1L_1\pr+a_2 L_2\pr+a_3\pr L_3\pr+a_3\prpr L_3\prpr+a_4\pr L_4\pr+a_4\prpr L_4\prpr+
eX\pr+f\pr Y\pr+f\prpr Y\prpr+g Z\pr)
$$
so that there is a covering of orbifolds
\begin{equation}\label{thirdlifting}
\cd(a_1,a_2;a_3,a_3;a_4,a_4; e;f,f;g)
\ok 
\cc_2(2a_1,2a_2,a_3,a_4;2e,f,g)
\end{equation}
In particular there is a covering
\begin{equation}\label{thirdlifting2}
\cd(1,1;2,2;2,2;1;2,2;1)
\ok
\cc_2(2,2,2,2;2,2)
\end{equation}
The locus of the orbifold $\cd(1,1;2,2;2,2;1;2,2;1)$ 
consists of the lines $\{L_{3}\pr,L_{3}\prpr,L_{4}\pr,L_{4}\prpr,Y\pr,Y\prpr\}$. 
Since the lifting $\gamma_2\1(L_3\cap L_4\cap Y)$  consist of four triple points, these lines 
forms an arrangement of 6 lines with 4 triple points, which is the complete quadrilateral.
In other words, one has an isomorphism $\cd(1,1;2,2;2,2;1;2,2;1)\simeq \cc_2(2,2,2,2;2,2,1))$.
The orbifolds 
$\cc_2(2,2,2,2;2e,2f,g;b_4,\dots,b_n)$ can be lifted recursively by these coverings 
$\cc_2(2,2,2,2;2,2,1)\ok\cc_2(2,2,2,2;2,2,1)$. The lifted orbifolds provide many examples of orbifolds 
$(\proj,D)$ uniformized by $\pball$. In Figure 9, we have shown an
orbifold obtained this way, namely the orbifold $\cd(1,1;2,2;2,2;1;2,2;2)$. Now one can 
take another set of three lines that do not meet at a triple point and lift this orbifold to 
the corresponding covering $\proj\ok\proj$. 

\begin{figure}
\label{unifbyb1b12}
\begin{center} 
\resizebox{0.7\hsize}{!}{\includegraphics*{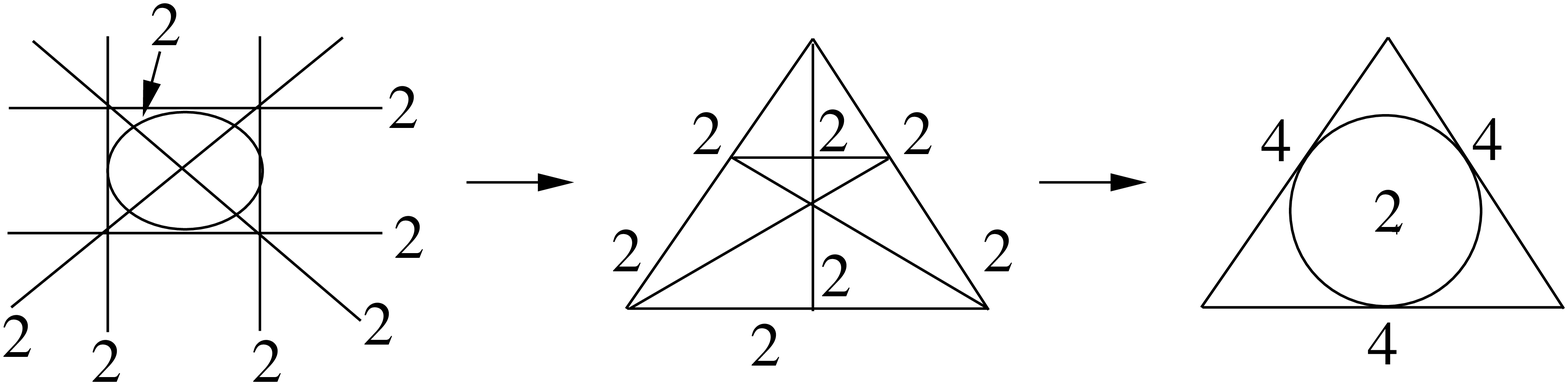}}  
\end{center}
\caption{Coverings $\cd(1,1;2,2;2,2;1;2,2;2)\ok\cc_2(2,2,2,2;2,2,2)\ok \ca(2;4,4,4)$} 
\end{figure} 
\par
Consider the covering  $\cc_4(2,2,2,2)\ok \ca(2;4,4,4)$.
The curves $Q_{4,1}$, $Q_{4,2}$, $Q_{4,3}$, $Q_{4,4}$ are smooth quadrics, 
such that $Q_{4,i}\cap Q_{4,j}$ consists of exactly two points, and 
$Q_{4,i}\cap Q_{4,j}\cap Q_{4,k}=\emptyset$ for $1\leq i\neq j\neq k\leq 4$.
Explicitly, these quadrics can be  given by the equations
$\alpha\pr x^2\underline{+}\beta\pr y^2\underline{+}\gamma\pr z^2=0$.
Since the orbifold $\ca(2;4,4,4)$ is uniformized by $\pball$, so is 
the orbifold $\cc_4(2,2,2,2)$.
\par
The orbifolds $\ca(2;b_1,\dots, b_n)$ can be lifted  to the uniformizations 
by K3 surfaces $X_i$ of the orbifolds given in the following lemma, yielding 
infinitely many examples of orbifolds $(X_i, D)$ uniformized by
$\pball$, with $X_i$ ($i\in\{1,2,3\}$) being a K3 surface.
\begin{lemma}
Let $T_i$ ($1\leq i\leq 6$) be six lines in $\proj$ in general position.
Then the orbifolds 
$$
\begin{array}{c}
\ce_1:=(\proj, 6T_1+6T_2+6T_3+2T_4),\\ 
\ce_2:=(\proj,4T_1+4T_2+4T_3+4T_4), \\ 
\ce_3:=(\proj, 2T_1+2T_2+2T_3+2T_4+2T_5+2T_6)
\end{array}
$$ 
are uniformized by K3 surfaces. 
\end{lemma}

\begin{figure}[h]\label{linek3s} 
\begin{center}
\resizebox{0.17\hsize}{!}{\includegraphics*{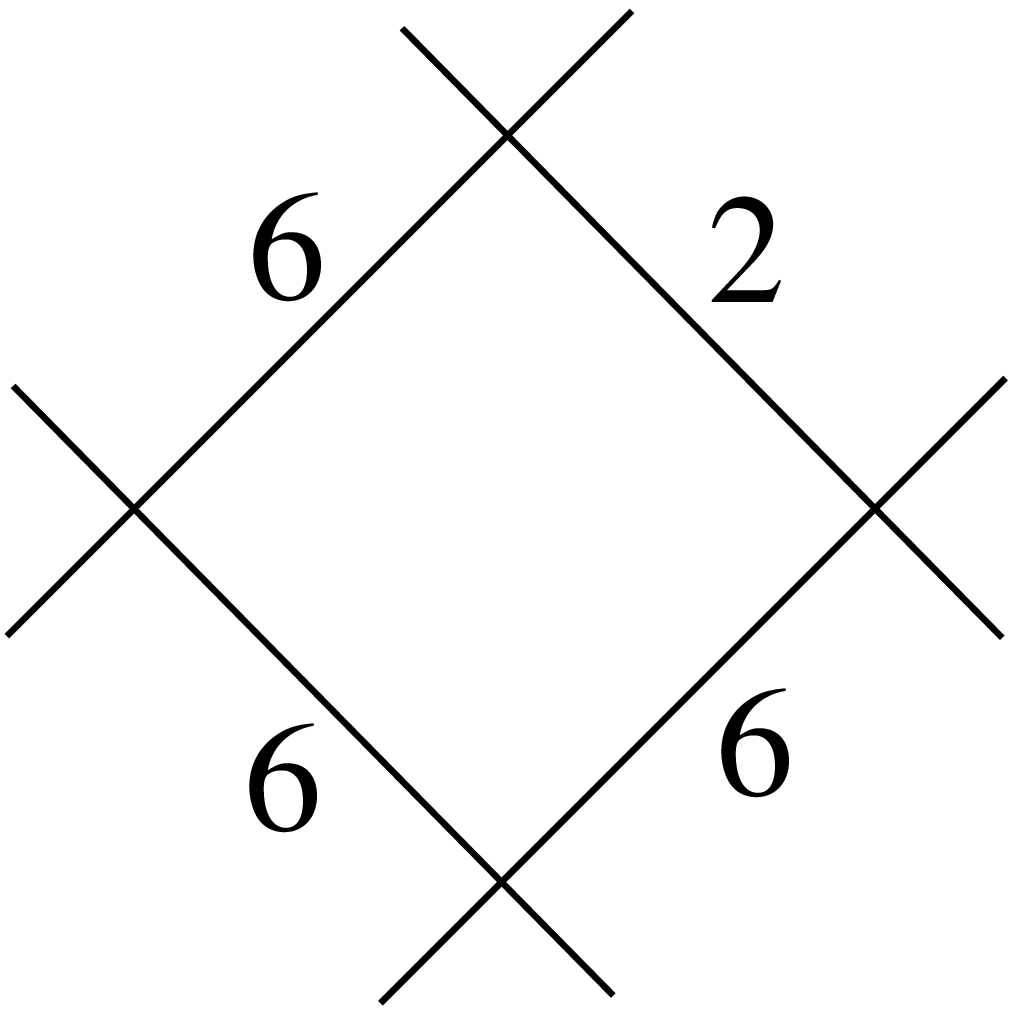}}  
\hspace{4mm} 
\resizebox{0.17\hsize}{!}{\includegraphics*{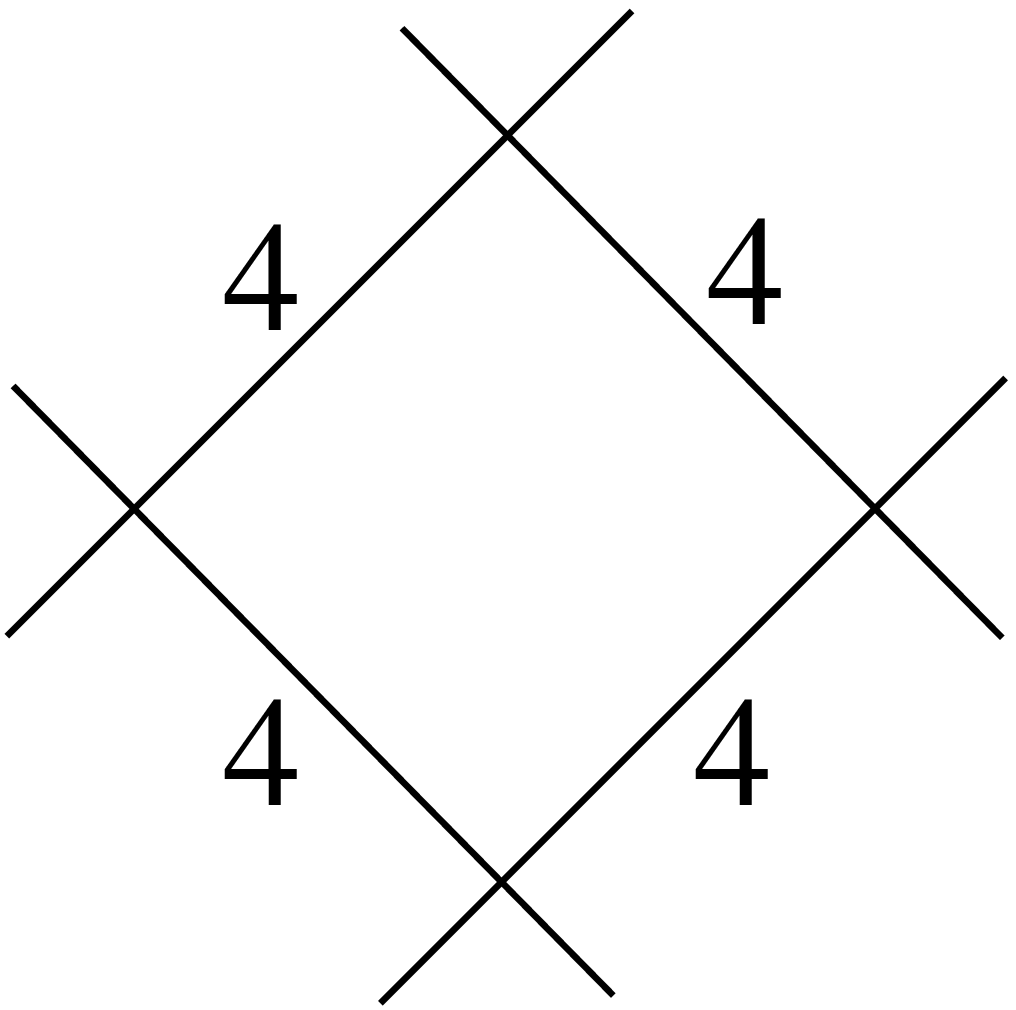}}  
\hspace{4mm} 
\resizebox{0.17\hsize}{!}{\includegraphics*{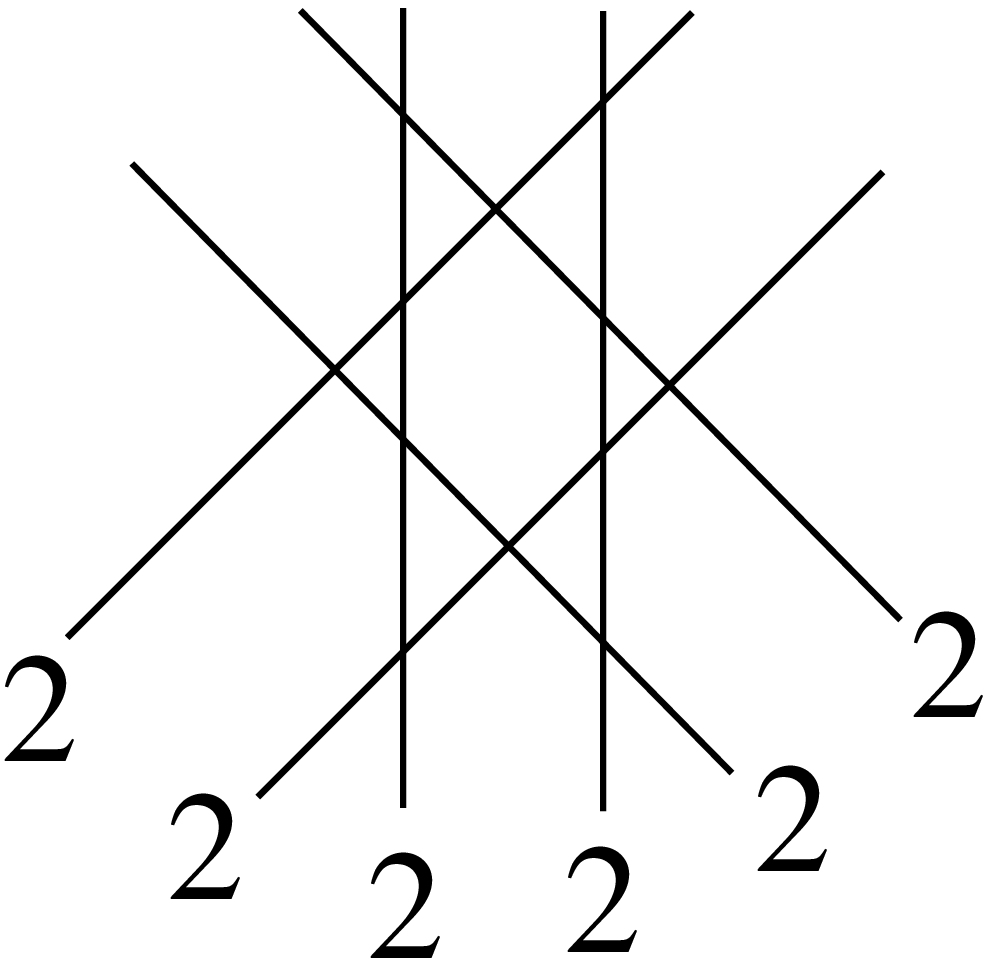}} 
\hspace{4mm} 
\end{center} 
\caption{Orbifolds uniformized by K3 surfaces} 
\end{figure} 

\begin{proof}
Lifting $\ce_1$ to the uniformization of $(\proj, 6T_1+6T_2+6T_3)$
yields an orbifold $(\proj,2T_{4,6})$ where $T_{4,6}$ is a smooth sextic.
This latter orbifold clearly admits a universal uniformization by a K3 surface.
Similarly, lifting $\ce_2$ to the uniformization of $(\proj, 4T_1+4T_2+4T_3)$
yields an orbifold $(\proj,4T_{4,4})$ where $T_{4,4}$ is a smooth quartic,
which clearly admits a universal uniformization by a K3 surface.
 The case of the orbifold $\ce_3$ is well known, see~ e.g. \cite{hunt}.

\end{proof}

\medskip\noindent\textbf{Proof of Theorem 2}
The orbifold $(\proj,2Q_m)$ is a lifting 
of the orbifold $\ca(2;m,m,m)$ to the uniformization 
of $(\proj,mT_1+mT_2+mT_3)$. For $i\in\{1,2,3\}$, let 
$p_i:=Q\cap T_i$. Then the orbifold $(\proj,mT_1+mT_2+mT_3)$ restricts
to $Q$ as the orbifold 
$\cs(m,m,m):=(Q,mp_1+mp_2+mp_3)$, provided that $m$ is an odd 
integer. The restriction of $\phi_m$ to $Q$ shows that 
there is a uniformization $Q_m\ok \cs(m,m,m)$ 
by an abelian map of degree $m^2$, {and consequently there is a uniformization
of degree $m^4$}
$$
\zeta :Q_m\times Q_m\ok \cs(m,m,m) \times \cs(m,m,m)
$$
On the other hand, with the notations of Section~\ref{apoconf} there is a covering
$$
\psi:\cs(m,m,m)\times \cs(m,m,m)\ok \ca(2;m,m,m)
$$
Taking the above uniformization of $\cs(m,m,m)$, we see that there 
is a Galois covering $\psi\comp \zeta:Q_m\times Q_m\ok \ca(2;m,m,m)$.
The map $\psi\comp \zeta$ is of degree $2m^4$. Let $G$ be the corresponding Galois group.
$$
\begin{array}{ccc}
(\proj,2Q_m)&\stackrel{\sigma}{\longleftarrow\!\!\!-\!\!-}&Q_m\times Q_m \\
&&\\
\phi_m\downarrow&&\downarrow\zeta\\
&&\\
\ca(2;m,m,m)&\stackrel{\psi}{\longleftarrow\!\!\!-\!\!-}&\cs(m,m,m)\times \cs(m,m,m)\\
\end{array}
$$
We want to show that the uniformization 
$Q_m\times Q_m\ok\ca(2;m,m,m)$ 
lifts to a uniformization of 
$(\proj,2Q_m)$. This is equivalent to showing that 
$$
H:=(\psi\comp\zeta )_*(\pi_1(Q_m)\times \pi_1(Q_m))\vartriangleleft 
K:=(\phi_m)_*(\orbfg(\proj, 2Q_m))
$$
The quotient of $\orbfg(\ca(2;m,m,m))$ by the normal subgroup generated 
by the meridians of $Q$ gives the group 
$$
\orbfg(\proj, mT_1+mT_2+mT_3)\simeq \Z/(m)\oplus \Z/(m)
$$
Since this is the Galois group of $\phi_m$, 
the group $K$ is the normal subgroup of $\orbfg(\ca(2;m,m,m))$ 
generated by the meridians of $Q$. Note that since 
$\orbfg(\ca(2;m,m,m))/K$ is abelian, $K$ should contain the commutators
$[\t_i,\t_j]$, for all meridians $\t_i$ of $T_i$ and $\t_j$ of $T_j$.
On the other hand, the quotient of $\orbfg(\ca(2;m,m,m))$ by the normal 
subgroup  generated by the meridians of $T_1$, $T_2$ and $T_3$ gives the group
$\orbfg(\proj, 2Q)\simeq \Z/(2)$,
which shows that  $\psi_*(\orbfg(S(m,m,m)\times S(m,m,m)))$ is the normal 
subgroup of $\ca(2;m,m,m)$ generated by the meridians of $T_1$, $T_2$ and 
$T_3$. Since 
$$
\zeta: Q_m\times Q_m\ok S(m,m,m)\times S(m,m,m)
$$
is the maximal abelian covering, the group $H$ is the normal subgroup of 
$\ca(2;m,m,m)$ generated by the commutators of the meridians of 
$T_1$, $T_2$ and $T_3$. This proves that $H$ is a normal subgroup of $K$. \QED

\medskip
Double covers of $\proj$ branched along the curves $Q_m$ 
were also studied in~\cite{persson}.
For a discussion of the groups $\fg{Q_m}$ see Section 7.
\section{Coverings of the orbifolds $\ca$ uniformized by $\ball$}
\noindent
\textbf{Coverings of $\ca(3;6,3,3)$.}\\
There is a covering $\phi_3:\cc_3(3;2,1,1)\ok\ca(3;3,3,6)$.
More explicitly, $\cc_3(3,2,1,1)$ is the orbifold 
$(\proj, 3Q_m+2X)$, where $Q_m$ is the nine-cuspidal 
sextic, and $X$ is a line passing through three of its cusps.

\medskip\noindent
\textbf{Coverings of $\ca(3;6,6,2)$.}\\
There is a covering $\phi_2:\cc_2(3,3,3,3;3,3,1)\ok \ca(3;6,6,2)$.
More explicitly, one has, as in Section~\ref{sectioncoverings},
$$
\cc_2(3,3,3,3;3,3,1):=(\proj, 3L_1+3L_2+3L_3+3L_4+3X+3Y), 
$$ 
where the set of lines $\{L_1,L_2,L_3,L_4,X,Y\}$ forms a
complete quadrilateral. Keeping the notations of Section~\ref{sectioncoverings},
consider the uniformization map
$\gamma_3:\proj\ok (\proj, 3X+3Y+3L_1)$.
The lifting $\gamma_3\1(L_2)$ consists of three lines 
$L_{2,1}$,  $L_{2,2}$ and $L_{2,3}$ which meet at the point 
$\gamma_3\1(L_{1,2}\cap L_{2,2} \cap Y)$. Similarly, the lifting 
$\gamma_3\1(L_{4,2})$ consists of three lines 
$L_{4,1}$, $L_{4,2}$, $L_{4,3}$ which meet at the point 
$\gamma_3\1(L_{4,2}\cap L_{1,2}\cap X)$. 
The lifting $\gamma_3\1(L_{3,2})$ is a smooth cubic $L_{3,1}$, and
it is readily seen by local considerations that 
for $i\in\{1,2,3\}$ 
the lines $L_{2,i}$, $L_{4,i}$ are tangent to $L_{3,1}$ with 
multiplicity 3. The points $\gamma_3\1(L_{2,2}\cap L_{4,2})$ lifts as the 
nine points of intersection $L_{2,i}\cap L_{2,j}$ ($i,j\in\{1,2,3\}$).
This shows that the lift of $\cc_2(3,3,3,3;3,3,1)$ along $\gamma_3$ is the orbifold
$$
(\proj, 3L_{3,1}+3L_{2,1}+3L_{2,2}+3L_{2,3}+3L_{4,1}+3L_{4,2}+3L_{4,3})
$$
with 2 points of type $x^3=y^3$, nine points of type $x^2=y^2$, and 
6 points of type $x^6=y^2$.
\par
Lifting this orbifold once more to the uniformization of 
$(\proj, 3L_{2,1}+3L_{2,2}+3L_{4,1})$ 
yields an arrangement of a smooth curve $L_{3,1,1}$ of degree 9 
with two smooth cubics $L_{4,2,1}$ and $L_{4,3,1}$ 
and three lines  $L_{2,3,1}$, $L_{2,3,2}$, $L_{2,3,3}$.
These lines meet at a point, and intersects the two cubics at $18$ 
nodes. The cubics are tangent to each other at 3 points with 
multiplicity 3. Cubics and lines altogether intersect the degree-9 curve 
at 27 distinct points of type $x^6=y^2$.  
\par
Returning again to the orbifold $\cc_2(3,3,3,3;3,3,1)$, consider the uniformization
map $\sigma_3:\proj\ok(\proj,3L_1+3L_2+3X)$. 
One has
$$
\sigma_3-lift(\cc(3,3,3,3;3,3,1))=(\proj,3K_1+3K_2+\dots +3K_9),
$$ 
where the lines $K_1,\dots,K_9$ forms a Ceva arrangement, which can be 
given by the equation
$(x^3-y^3)(y^3-z^3)(z^3-x^3)=0$.
Suppose that $K_1$, $K_2$, $K_3$ do not meet at a triple point.
Lifting this orbifold to the uniformization by $\proj$ 
of $(\proj,3K_1+3K_2+3K_3)$ yields an arrangement of nine lines 
with three smooth cubics. It is left to the reader to verify that this 
arrangement can be lifted once more to $\proj$ in two different ways.

\medskip
\noindent\textbf{Coverings of $\ca(3;3,4,4)$.}\\
This orbifold (and the orbifold $\ca(3;6,6,2)$) can be lifted to a K3 - uniformization  of 
 $\ca(3;3,2,2)$. An example of a ball quotient orbifold over a K3 surface was 
also given in~\cite{naruki}. Recall that $\ca(3;3,4,2)$
satisfy $c_1^2=2e$, but we don't know whether it admits a 
uniformization. In case it does, $\ca(3;3,4,4)$ lifts to this uniformization.

\medskip\noindent
\textbf{Coverings of $\ca(4;4,4,4)$.}\\
Lifting $\ca(4;4,4,4)$ by $\phi_4$ gives the 
orbifold $\cc_4(4,4,4,4)$ over four mutually tangent quadrics.
It can also be lifted to the uniformization of $\ca(4;2,2,2)$, 
which is a K3 surface. Note that by Theorem~\ref{firstlifting} the orbifold
$\ca(2;4,4,4)$ is uniformized by a product of two Riemann surfaces.
The orbifold  $\ca(4;4,4,4)$ can also be lifted to this product. 
\begin{figure}

\begin{center} 
\resizebox{0.5\hsize}{!}{\includegraphics*{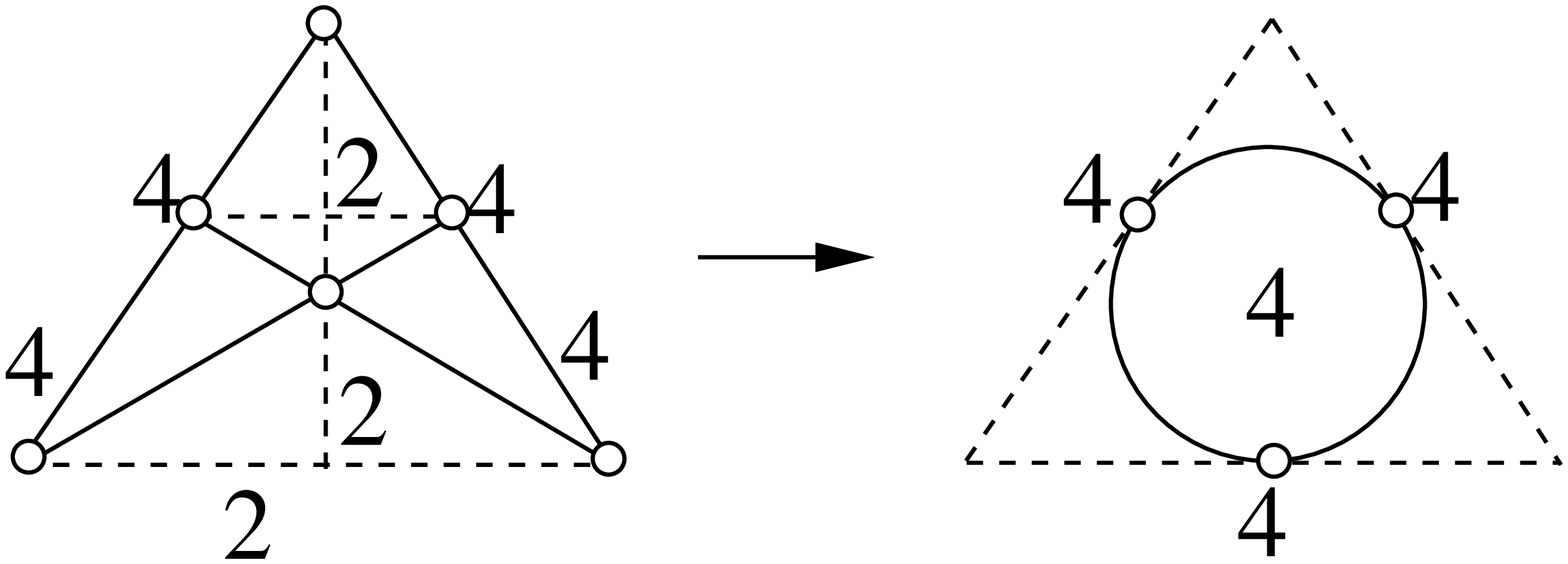}}  
\end{center}
\caption{The covering $\cc_2(4,4,4,4;2,2,2)\ok\ca(4;4,4,4)$} 
\end{figure} 

\medskip\noindent\textbf{Proof of Theorem 1}
There is an orbifold covering 
$$\mu_1:\cc_2(4,4,4,4;2,2,2)\ok\ca(4;4,4,4),$$
where $\mu_1$ is a bicyclic covering of degree 4, 
branched along the dashed lines in the locus of $\ca(4,4,4,4)$ 
(see Figure 11). 
The lattices corresponding to 
both of the orbifolds $\ca(4;4,4,4)$ and $\co_1:=\cc_2(4,4,4,4;2,2,2)$ are known to be 
arithmetic, see~\cite{holzapfel3} and~\cite{deligne2}.
In the locus of the orbifold  $\co_1$, take three dashed lines and mark the remaining
lines with red and blue as in Figure~12.  
Consider the degree-4 bicyclic covering $\mu_2:\proj\ok\proj$ branched along the dashed lines.

\begin{figure}[h]

\begin{center} 
\resizebox{0.5\hsize}{!}{\includegraphics*{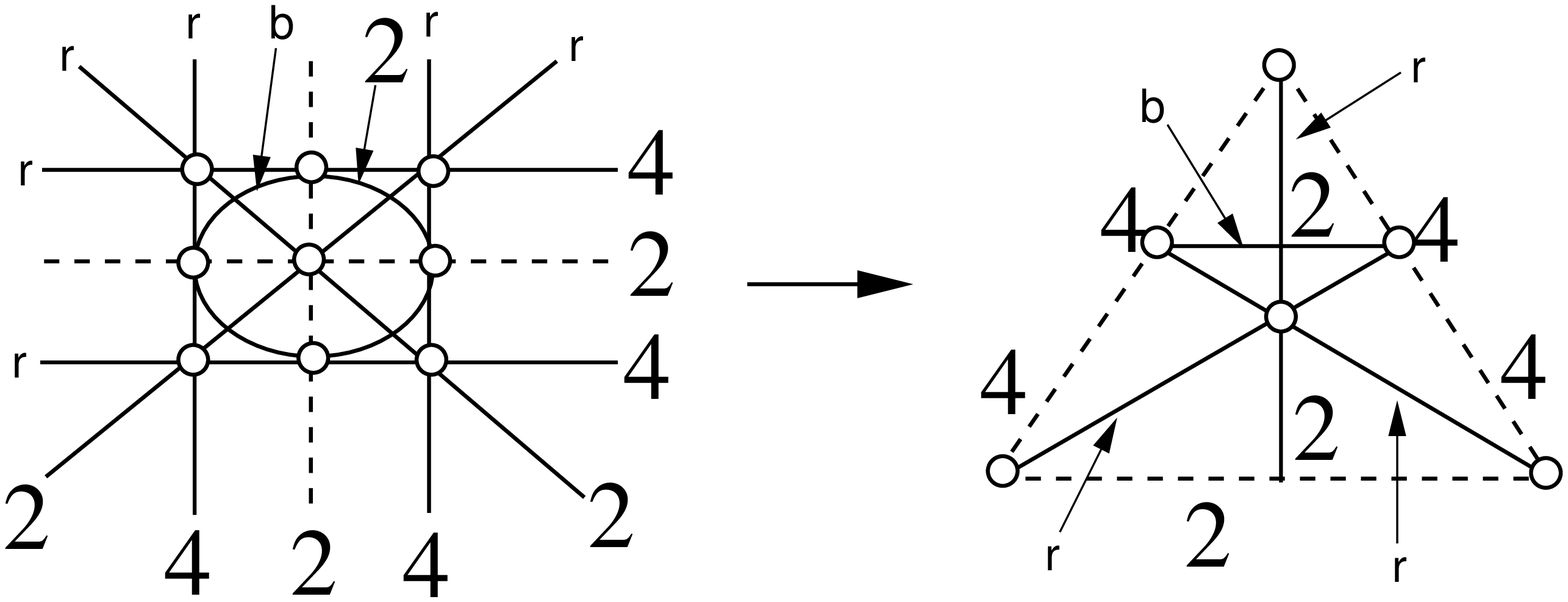}}  
\end{center}
\caption{The covering $\co_2\ok \co_1:=\cc_2(4,4,4,4;2,2,2)$} 
\end{figure} 

Since the orbifold $\co_1$ is uniformized by $\ball$, so is the orbifold 
$\co_2:=\mu_2-lift\,(\co_1)$. Since the red lines pass through the intersection points of the 
dashed lines, the red lines 
 will be lifted as 6 lines forming a complete quadrilateral. Since the blue line 
intersects the dashed lines at three distinct points, its lifting will be a smooth quadric. 
Two  dashed lines with weight 4 also lifts as two lines, and their weight becomes 2 
(see Figure 12). Now redraw the locus of the orbifold and mark the curves in this locus 
as in Figure 13.
\begin{figure}[h]

\begin{center} 
\resizebox{0.25\hsize}{!}{\includegraphics*{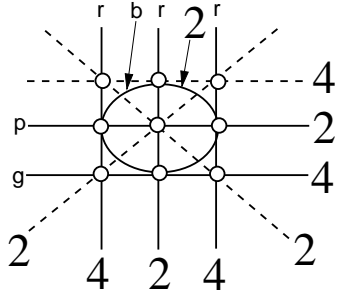}}  
\end{center}
\caption{New marking of the orbifold $\co_2$} 
\end{figure} 

Consider the degree-4 bicyclic covering $\mu_3:\proj\ok\proj$ branched along the dashed lines.
Since the orbifold $\co_1$ is uniformized by $\ball$, so is the orbifold 
$\co_3:=\mu_3-lift\,(\co_2)$. Let us describe the locus of $\co_3$. Since the red lines pass through
the intersection points of the dashed lines, they will lift as 6 lines, forming a complete 
quadrilateral. Denote these lines by $R_1,\dots,R_6$. The pink line also pass through 
an intersection point of the dashed lines, so it will lift as two lines $P_1$ and $P_2$. 
Green line intersects the dashed lines at three distinct points, so its lifting will be a 
smooth quadric $G$. The blue quadric will be lifted as a quartic curve $B$ 
(which is in fact irreducible with two nodes). The dashed line with weight four 
will also lift as a line $D$ with weight 2. Hence, $\co_3$ is the orbifold
$$
(\proj, 4R_1+4R_2+4R_3+4R_4+2R_5+2R_6+2P_1+2P_2+4G+2B+2D)
$$
Consider the complete quadrilateral formed by the lines $R_1,\dots,R_6$.
Suppose that the lines $R_1, R_3, R_5$ meets at a point. Then the lines $R_2,R_4,R_6$
intersect each other at three distinct points. Let $\mu_4$ be the degree-4 bicyclic covering 
$\mu_4:\proj\ok\proj$ branched along $R_2,R_4,R_6$. Consider the orbifold 
$\co_4:=\mu_4-lift\,(\co_3)$. Since the lines $R_1,R_3,R_5$ pass through the 
intersection points of $R_2,R_4,R_6$, they fill be lifted as 6 lines forming a complete 
quadrilateral, so that the locus of $\co_4$ will contain a complete quadrilateral. 
Now one can recursively apply this proces to get an infinite series $\co_r$ 
of orbifolds uniformized by $\ball$.
Let $\mu_r:\co_{r}\ok\co_{r-1}$ be the $r+1$th covering in the recursion. 
Consider the blue quadric in the  locus of $\co_2$.
Mark a curve in the locus of  $\co_r$ with blue if $\mu_r(C)$ is blue. 
During the recursion, one may assume that 
$\mu_r$ is never branched along a blue curve. Such a covering has the effect of multiplying
the total degree of blue curves by two. Hence, the degree of the locus of $\co_r$ is 
$\geq 2^r$, which shows that among  $\co_r$ there are infinitely many distinct orbifolds.
\QED

\medskip \noindent\textbf{Remarks.} (1) For $r\geq2$, the uniformizations $\ball\ok\co_r$ take place outside 
of finitely many cusp points. Indeed, the orbifold $\co_2$ has cusp points only 
and locally an orbifold-covering of a cusp point is also a cusp point. But since  
$\co_r\ok\co_{r-1}$ is an orbifold covering,  singular points of $\co_r$ are locally coverings
of the singular points of $\co_{r-1}$. \\
(2) One may apply the recursion described above in many different ways, for example taking three 
lines in the locus of $\co_2$, which are tangent to the blue quadric and lifting $\co_2$ 
to the bicyclic covering $\proj\ok\proj$ branched along these lines gives another  infinite 
series of ball-quotient orbifolds.\\
(3) The orbifold $\ca(4;4,4,4)$ is invariant under a $\Sigma_3$-action on $\proj$, and
the quotient orbifold $\cq:=\ca(4;4,4,4)/\Sigma_3$ (with a singular base space $\proj/\Sigma_3$)
is the ``queen'' of the orbifolds constructed above.  Its locus consists of two rational curves,
one is the image of the quadric, and the other  is the image of the tangent lines.

\section{Fundamental groups}\label{section:7}
A presentation of the group $\fg{\ca_n}$ is known, see~\cite{amram} and 
see~\cite{paolo}, \cite{paolo2}, \cite{lambro} for generalizations.
\begin{theorem}\label{thepresentation}
One has the presentation
\begin{equation}\label{group1}
\fg{\ca_n} \simeq \Biggl\langle 
\begin{array}{l}\t_1,\dots,\t_n,\\
 \q_1,\dots,\q_n 
\end{array}\; \left | 
\begin{array}{ll}
 \q_i = \t_i \q_{i-1} \t_i^{-1},\  2 \leq i \leq n \\
(\q_i \t_i)^2 = (\t_i \q_i)^2,\   1 \leq i \leq n \\
\left[\q_i^{-1} \t_i \q_i, \t_j \right] = 1,\  1 \leq i<j\leq n\\
\t_n \cdots \t_1 \q_1^2 = 1 
\end{array} 
\right . \Biggr\rangle 
\end{equation}
where $\q_i$ are meridians of $Q$ and $\t_i$ is a meridian of 
$T_i$ for $1\leq i\leq n$. 
\end{theorem}
Note that $\fg{\ca_n}$ is the braid group on two strands of the punctured sphere
$Q\backslash\{p_1,\dots,p_n\}$. For small $n$, its presentation can be simplified: 
\begin{corollary}\label{simplified}\cite{amram}
(i)  The group $\fg{\ca_1}$ is abelian.\\
(ii) \cite{dethloff} The group $\fg{\ca_2}$ admits the presentation 
$\fg{\ca_2}\simeq \langle \t,\q\,|\, (\t\q)^2=(\q\t)^2 \rangle$,
where $\q$ is a meridian of $Q$ and $\t$ is a meridian of $T_1$. A meridian of
$T_2$ is given by $\q^{-2}\t\1$.\\
(iii)  \cite{degtyarev} The group $\fg{\ca_3}$ admits the presentation 
$\fg{\ca_3}\simeq 
\langle \t,\s,\q\,|\, (\t\q)^2=(\q\t)^2,\; (\s\q)^2=(\q\s)^2,\;[\s,\t]=1\rangle$,
where $\s$, $\t$ are meridians of $T_1$ and $T_3$ respectively, 
and $\q$ is a meridian of $Q$. A meridian of $T_2$ is given by 
$(\q\t\q\s)\1$.
\end{corollary}
\begin{corollary}\label{modular}
One has the presentation
$$
\orbfg(\ca(a;b_1,b_2,b_3))\simeq\left\sol \q,\t,\s\,
\begin{array}{|c}
(\t \q)^2=(\q\t)^2,\quad (\s \q)^2=(\q\s)^2,\\
\left[\s,\t\right]=\q^a=\t^{b_1}=\s^{b_2}=(\q\t\q\s)^{b_3}=1
\end{array}\right\sag
$$
\end{corollary}
\noindent
Fundamental groups of almost all the curves or arrangements 
appearing in this article is a subgroup of this group. 
In particular, there is an exact sequence
$$
0\ok\fg{Q_m}\ok\orbfg(\ca(\infty;m,m,m))\ok\Z/(m)\oplus\Z/(m)\ok 0
$$
In~\cite{cogolludo}, a presentation of 
$\fg{Q_m}$ was computed from this exact sequence 
by using the Reidemeister-Schreier algorithm.
A presentation of the group $\fg{Q_3}$ 
was found by Zariski~\cite{zariski2}, see also~\cite{kaneko}.

\section{Final remarks}
Let $C\subset \proj$ be an irreducible curve with $\kappa$ simple cusps $\nu$ 
nodes and no other singularities. 
The Euler number of the orbifold $(\proj, bC)$ ($b\in\{2,3,4,5,6\}$) is
$$
e(\proj,bC)=e(\proj\moins C)+\frac{e(C\moins\mbox{Sing}(C))}{b}+
\frac{3\kappa}{2}\left[\frac{1}{b}-\frac{1}{6}\right]^2
$$
Considering $C$ as a subset of $\proj$, one has  $e(C)=-d^2+3d+2\kappa+\nu$,
so that $e(\proj\moins C)=3+d^2-3d-2\kappa-\nu$. 
Setting $e(C\moins\mbox{Sing}(C))=e(C)-\kappa-\nu$ gives
$$
e(\proj,mC)=3+d^2-3d-2\kappa-\nu+\frac{-d^2+3d+\kappa}{b}+
\frac{\nu}{b^2}+\frac{3\kappa}{2}\left[\frac{1}{b}-\frac{1}{6}\right]^2
$$
On the other hand, the first Chern number of this orbifold is
$$ 
c_1^2(\proj,bC)=\left[-3+d\left(1-\frac{1}{b}\right)\right]^2
$$    
\par
Let $g$ be the genus of $C$. It is easy to show that there exists infinitely 
many five-tuples $(d,\kappa,\nu,b,g)$ with $(3e-c_1^2)(\proj,bC)=0$.
Some examples are given in the table below. 
The first curve in the table is the 9-cuspidal sextic. 
We don't know whether the other curves exists.

{\small
$$
\begin{array}{l|lllll||}
&d& \kappa & \nu& b& g \\ \hline
 \textbf{1}&6&  9&  0&  2&  1\\
  \textbf{2}&7&  11& 3&  2&  1\\
  \textbf{3}&8&  17&  0&  2&  4\\
  \textbf{4}&9&  18&  7&  2&  3\\
  \textbf{5}&9&  18&  8&  4&  2\\
  \textbf{6}&10&  23&  8&  2&  5\\
  \textbf{7}&11&  32&  3&  2&  10\\
  \textbf{8}&11&  23&  19&  2&  3\\
  \textbf{9}&12&  36&  8&  2&  11\\
  \textbf{10}&12&  36&  10&  4&  9
\end{array}\begin{array}{l|lllll||}
&d& \kappa & \nu& b& g \\ \hline
 \textbf{11}&12&  27&  24&  2&  4\\
 \textbf{12}& 13&  44&  7&  2&  15\\
 \textbf{13}&13&  35&  23&  2&  8\\
 \textbf{14}&13&  26&  39&  2&  1\\
 \textbf{15}&14&  56&  0&  2&  22\\
 \textbf{16}&14&  47&  16&  2&  15\\
 \textbf{17}&14&  38&  32&  2&  8\\
 \textbf{18}&14&  29&  48&  2&  1\\
 \textbf{19}&15&  63&  3&  2&  25\\
 \textbf{20}&15&  64&  3&  3&  24
\end{array}\begin{array}{l|lllll}
&d& \kappa & \nu& b& g \\ \hline
 \textbf{21}&15&  54&  19&  2&  18\\
 \textbf{22}&15&  45&  35&  2&  11\\
 \textbf{23}&15&  36&  51&  2&  4\\
 \textbf{24}&15&  40&  51&  6&  0\\
 \textbf{25}&16&  74&  0&  2&  31\\
 \textbf{26}&16&  65&  16&  2&  24\\
 \textbf{27}&16&  56&  32&  2&  17\\
 \textbf{28}&16&  47&  48&  2&  10\\
 \textbf{29}&16&  38&  64&  2&  3\\
 \textbf{30}&17&  80&  7&  2&  33
\end{array}
$$
}

\bigskip\noindent
{\tiny \textbf{Acknowledgements.}
This work was partially achieved during the author's post-doctorate stay in
the Emmy Noether Research Institute for Mathematics
(center of the Minerva Foundation in Germany), and was partially supported by 
the Excellency Center ``Group Theoretic
Methods in the Study of Algebraic Varieties'' of the Israel Science Foundation, and by 
EAGER (EU network, HPRN-CT-2009-00099). The author is grateful to Meirav Amram,
Prof. Mina Teicher and to the Emmy Noether Research Institute for their hospitality.  
The author also expresses his gratitude to the Referee of Math. Annalen for his remarks.}


\noindent
Galatasaray University, Department of Mathematics, 80840 
Ortak{\"o}y /Istanbul, TURKEY \\ email~:~\texttt{muludag@gsu.edu.tr}

\end{document}